\documentclass[11pt]{article}

	\newcommand{\mainTitle}{Use of the generating function to generalize the sum formula for quadruple zeta values}

	\newcommand{\authorName}{MACHIDE, Tomoya}
	\newcommand{\organizationNameFst}{National Institute of Informatics}
	\newcommand{\placeAddressFst}{2-1-2 Hitotsubashi, Chiyoda-ku, Tokyo 101-8430, Japan}
	\newcommand{\emailAddressFst}{machide@nii.ac.jp}
	\newcommand{\organizationNameScd}{JST, ERATO, Kawarabayashi Large Graph Project}
	\newcommand{\departmentNameScd}{Japan} 
	
	\newcommand{\MSCname}{11M32(Primary), 16S34,20C07(Secondary)} 
	\newcommand{\keyWord}{multiple zeta value, sum formula, generating function, group ring of general linear group}

\usepackage{graphicx}
\usepackage{geometry}               
\usepackage{ascmac}
\usepackage{amsmath}
\usepackage{amssymb}
\usepackage{amsthm}
\usepackage[colorlinks=true]{hyperref}
\usepackage{arydshln}
\usepackage{multicol}
	
	\DeclareMathOperator*{\OPlus}{\bigoplus}
	\newcommand{\nbk}[3]{#1#3#2}		
	\newcommand{\bgbk}[3]{\bigl{#1}#3\bigr{#2}}	
	\newcommand{\Bgbk}[3]{\Bigl{#1}#3\Bigr{#2}}			
	\newcommand{\bggbk}[3]{\biggl{#1}#3\biggr{#2}}			
	\newcommand{\Bggbk}[3]{\Biggl{#1}#3\Biggr{#2}}
	\newcommand{\autobk}[3]{\left#1#3\right#2}
	\newcommand{\nbkD}[5]{#1#2#5#3#4}		
	\newcommand{\bgbkD}[5]{\bigl{#1}\bigl{#2}#5\bigr{#3}\bigr{#4}}	
	\newcommand{\BgbkD}[5]{\Bigl{#1}\Bigl{#2}#5\Bigr{#3}\Bigr{#4}}	
	\newcommand{\bggbkD}[5]{\biggl{#1}\biggl{#2}#5\biggr{#3}\biggr{#4}}	
	\newcommand{\BggbkD}[5]{\Biggl{#1}\Biggl{#2}#5\Biggr{#3}\Biggr{#4}}	
	\newcommand{\autobkD}[5]{\left#1\left#2#5\right#3\right#4}	
	\newcommand{\mcbk}[4][?]{\ifx n#1\nbk{#2}{#3}{#4}\else\ifx b#1\bgbk{#2}{#3}{#4}\else\ifx B#1\Bgbk{#2}{#3}{#4}\else\ifx g#1\bggbk{#2}{#3}{#4}\else\ifx G#1\Bggbk{#2}{#3}{#4}\else\ifx a#1\autobk{#2}{#3}{#4}\else\ifx !#1{#4}\else#4\fi\fi\fi\fi\fi\fi\fi}
	\newcommand{\mcbkD}[4][?]{\ifx n#1\nbkD{#2}{#2}{#3}{#3}{#4}\else\ifx b#1\bgbkD{#2}{#2}{#3}{#3}{#4}\else\ifx B#1\BgbkD{#2}{#2}{#3}{#3}{#4}\else\ifx g#1\bggbkD{#2}{#2}{#3}{#3}{#4}\else\ifx G#1\BggbkD{#2}{#2}{#3}{#3}{#4}\else\ifx a#1\autobkD{#2}{#2}{#3}{#3}{#4}\else\ifx !#1{#4}\else#4\fi\fi\fi\fi\fi\fi\fi}
	\newcommand{\nsgsb}[1]{#1}		
	\newcommand{\bgsgsb}[1]{\big{#1}}	
	\newcommand{\Bgsgsb}[1]{\Big{#1}}			
	\newcommand{\bggsgsb}[1]{\bigg{#1}}			
	\newcommand{\Bggsgsb}[1]{\Bigg{#1}}
	\newcommand{\mcsgsb}[2][?]{\ifx n#1\nsgsb{#2}\else\ifx b#1\bgsgsb{#2}\else\ifx B#1\Bgsgsb{#2}\else\ifx g#1\bggsgsb{#2}\else\ifx G#1\Bggsgsb{#2}\else#2\fi\fi\fi\fi\fi}
	\newcommand{\myEqSpace}{\,} 	\newlength{\myEqSpaceLen} 	
	\newcommand{\mLt}[1]{\widetilde{#1}}
	\newcommand{\pLt}[1]{\check{#1}}

	\newcommand{\bkR}[2][n]{\mcbk[#1]{(}{)}{#2}}						
	\newcommand{\bkS}[2][n]{\mcbk[#1]{[}{]}{#2}}						\newcommand{\bkSd}[2][n]{\mcbkD[#1]{[}{]}{#2}}
	\newcommand{\bkB}[2][n]{\mcbk[#1]{\{}{\}}{#2}}						
	\newcommand{\bkA}[2][n]{\mcbk[#1]{\langle}{\rangle}{#2}}				
	\newcommand{\bkAll}[4][n]{\mcbk[#1]{#2}{#3}{#4}}
		
	\newcommand{\nFc}[3][n]{#2\bkR[#1]{#3}}					
				
	\newcommand{\idFc}[4][n]{\id{#2}{#3}\bkR[#1]{#4}}			
	\newcommand{\pwFc}[4][n]{\pw{#2}{#3}\bkR[#1]{#4}}			
	\newcommand{\ipFc}[5][n]{\ip{#2}{#3}{#4}\bkR[#1]{#5}}		
		\newcommand{\Fc}{\nFc}		
			
	\newcommand{\nFcs}[3][n]{#2\bkS[#1]{#3}}

		\newcommand{\Fcs}{\nFcs}
		
	\newcommand{\gam}{\gamma} 
	\newcommand{\Gam}{\Gamma} 
	
	\newcommand{\ep}{\varepsilon}

	\newcommand{\sig}{\sigma}
	\newcommand{\Sig}{\Sigma}
	
	\newcommand{\Ome}{\Omega}

	\newcommand{\bfAlp}[1][s]{\ifx s#1{\boldsymbol\alpha}\else{\boldsymbol??none??}\fi}		\newcommand{\bfAA}{\bfAlp}
	\newcommand{\bfBeta}[1][s]{\ifx s#1{\boldsymbol\beta}\else{\boldsymbol??none??}\fi}		\newcommand{\bfBB}{\bfBeta}
	\newcommand{\bfDelta}[1][s]{\ifx s#1{\boldsymbol\delta}\else{\boldsymbol??none??}\fi}		
	\newcommand{\bfGam}[1][s]{\ifx s#1{\boldsymbol\gam}\else{\boldsymbol\Gam}\fi}			
	\newcommand{\bfA}[1][s]{\ifx s#1{\bf a}\else{\bf A}\fi}
	\newcommand{\bfB}[1][s]{\ifx s#1{\bf b}\else{\bf B}\fi}
	\newcommand{\bfE}[1][s]{\ifx s#1{\bf e}\else{\bf E}\fi}
	\newcommand{\bfH}[1][s]{\ifx s#1{\bf h}\else{\bf H}\fi}
	\newcommand{\bfI}[1][s]{\ifx s#1{\bf i}\else{\bf I}\fi}
	\newcommand{\bfJ}[1][s]{\ifx s#1{\bf j}\else{\bf J}\fi}
	\newcommand{\bfK}[1][s]{\ifx s#1{\bf k}\else{\bf K}\fi}
	\newcommand{\bfL}[1][s]{\ifx s#1{\bf l}\else{\bf L}\fi}
	\newcommand{\bfM}[1][s]{\ifx s#1{\bf m}\else{\bf M}\fi}
	\newcommand{\bfN}[1][s]{\ifx s#1{\bf n}\else{\bf N}\fi}
	\newcommand{\bfS}[1][s]{\ifx s#1{\bf s}\else{\bf S}\fi}
	\newcommand{\bfU}[1][s]{\ifx s#1{\bf u}\else{\bf U}\fi}
	\newcommand{\bfV}[1][s]{\ifx s#1{\bf v}\else{\bf V}\fi}
	\newcommand{\bfW}[1][s]{\ifx s#1{\bf w}\else{\bf W}\fi}
	\newcommand{\bfX}[1][s]{\ifx s#1{\bf x}\else{\bf X}\fi}
	\newcommand{\bfY}[1][s]{\ifx s#1{\bf y}\else{\bf Y}\fi}
	\newcommand{\bfZ}[1][s]{\ifx s#1{\bf z}\else{\bf Z}\fi}
	\newcommand{\bfEp}[1][s]{\ifx s#1{\boldsymbol\ep}\else{\boldsymbol{\mathcal{E}}}\fi}
	
	\newcommand{\mVert}[1][n]{{\,\mcsgsb[#1]{\vert}\,}}		
	
	\newcommand{\SetO}[2][n]{\bkB[#1]{#2}}
	\newcommand{\SetT}[3][n]{\bkB[#1]{#2\mVert#3}}
		\newcommand{\Set}{\SetO}
	\newcommand{\SpO}[2][n]{\bkA[#1]{#2}}

	\newcommand{\GpO}[2][n]{\bkA[#1]{#2}}
		\newcommand{\Gp}{\GpO}
	\newcommand{\setN}{\mathbb{N}}
	\newcommand{\setZ}{\mathbb{Z}} 
	\newcommand{\setQ}{\mathbb{Q}}	
	\newcommand{\setR}{\mathbb{R}}	
	
	\newcommand{\setF}[1][?]{\ifx ?#1\mathbb{F}\else\mathbb{F}_{#1}\fi}
	
	\newcommand{\matI}[1][?]{\ifx #1?I\else I_{#1}\fi}	
	
	\newcommand{\gpSym}[2][?]{S_{#2}}	
	\newcommand{\gpAlt}[2][?]{A_{#2}}		
	\newcommand{\gpKleinF}[1][?]{V}
		
	\newcommand{\gpGL}[3][?]{GL_{#2}(#3)}	
	\newcommand{\gpu}[1][?]{\ifx?#1e\else e_{#1}\fi}		
	\newcommand{\rgMat}[2]{M_{#1}(#2)}

	\newcommand{\vPack}[1][10]{\vspace{-#1pt}}
	
	\newcommand{\lnA}[1][]{&  &}
	\newcommand{\lnP}[1]{\myEqSpace#1\myEqSpace}
	\newcommand{\lnAP}[2][]{& #2 &}
	\newcommand{\lnAH}[1][\nonumber]{#1\\ & &}
	\newcommand{\lnAHP}[2][\nonumber]{#1 \\ & #2 &}
	\newcommand{\lnAHs}[2][\nonumber]{#1 \\ & &\hspace{#2pt}}
	
		\newcommand{\slnAH}[1][?]{\\}
		
	\newcommand{\refEq}[1]{(\ref{#1})}	
	\newcommand{\refEqA}[1]{(#1)}	
	\newcommand{\pcstSpForRefThm}{\;}		
	\newcommand{\refHL}[2]{#1\pcstSpForRefThm\ref{#2}}		
	\newcommand{\refHLm}[3][?]{\ifx?#1#2\pcstSpForRefThm#3\else#2#3\fi}
	\newcommand{\refThm}[2][?]{\ifx?#1\refHL{Theorem}{#2}\else\ifx s#1\refHL{Theorems}{#2}\else{[argument error]}\fi\fi}
	\newcommand{\refProp}[2][?]{\ifx?#1\refHL{Proposition}{#2}\else\ifx s#1\refHL{Propositions}{#2}\else{[argument error]}\fi\fi}
	\newcommand{\refLem}[2][?]{\ifx?#1\refHL{Lemma}{#2}\else\ifx s#1\refHL{Lemmas}{#2}\else{[argument error]}\fi\fi}
	\newcommand{\refCor}[2][?]{\ifx?#1\refHL{Corollary}{#2}\else\ifx s#1\refHL{Corollaries}{#2}\else{[argument error]}\fi\fi}
	\newcommand{\refDef}[2][?]{\ifx?#1\refHL{Definition}{#2}\else\ifx s#1\refHL{Definitions}{#2}\else{[argument error]}\fi\fi}
	\newcommand{\refRem}[2][?]{\ifx?#1\refHL{Remark}{#2}\else\ifx s#1\refHL{Remarks}{#2}\else{[argument error]}\fi\fi}
	\newcommand{\refTab}[2][?]{\ifx?#1\refHL{Table}{#2}\else\ifx s#1\refHL{Tables}{#2}\else{[argument error]}\fi\fi}
	\newcommand{\refSec}[2][?]{\ifx?#1\refHL{Section}{#2}\else\ifx s#1\refHL{Sections}{#2}\else{[argument error]}\fi\fi}
	\newcommand{\refApp}[2][?]{\ifx?#1\refHL{Appendix}{#2}\else\ifx s#1\refHL{Appendices}{#2}\else{[argument error]}\fi\fi}
		
	\newcommand{\refThmA}[2][?]{\ifx?#1\refHLm{Theorem}{#2}\else \refHLm{Theorems}{#2}\fi}
	\newcommand{\refPropA}[2][?]{\ifx?#1 \refHLm[#1]{Proposition}{#2}\else \refHLm{Propositions}{#2}\fi}
	\newcommand{\refLemA}[2][?]{\ifx?#1\refHLm{Lemma}{#2}\else \refHLm{Lemmas}{#2}\fi}
	\newcommand{\refCorA}[2][?]{\ifx?#1\refHLm{Corollary}{#2}\else \refHLm{Corollaries}{#2}\fi}
	\newcommand{\refDefA}[2][?]{\ifx?#1\refHLm{Definition}{#2}\else \refHLm{Definitionss}{#2}\fi}
	\newcommand{\refRemA}[2][?]{\ifx?#1\refHLm{Remark}{#2}\else \refHLm{Remarks}{#2}\fi}
	\newcommand{\refSecA}[2][?]{\ifx?#1\refHLm{Section}{#2}\else \refHLm{Sections}{#2}\fi}
	\newcommand{\refAppA}[2][?]{\ifx?#1\refHLm{Appendix}{#2}\else \refHLm{Appendices}{#2}\fi}
		\newcommand{\refSectA}{\refSecA}
	
	\newcommand{\frc}[3][?]{\ifx s#1#3/#2\else\ifx b#1(#3)/#2\else\ifx d#1\dfrac{#3}{#2}\else\ifx t#1\tfrac{#3}{#2}\else\frac{#3}{#2}\fi\fi\fi\fi}
	\newcommand{\mopI}[2][?]{\bkR[#1]{#2}^{-1}}
	
	\newcommand{\mopTP}[2][?]{\bkR[#1]{#2}^t}		

	\newcommand{\racFT}[3][?]{#2\vert\bkR[#1]{#3}}						\newcommand{\racFrrT}[3][?]{(#2\vert\bkR[#1]{#3})}	
												\newcommand{\racF}{\racFT}			\newcommand{\racFrr}{\racFrrT}
	\newcommand{\racFATh}[4][n]{\bkR[#1]{#2\vert#3}(#4)}				
												\newcommand{\racFA}{\racFATh}				
	\newcommand{\racFAsTh}[4][n]{\bkR[#1]{#2\vert#3}[#4]}				
												\newcommand{\racFAs}{\racFAsTh}

	\newcommand{\pw}[3][?]{\ifx!#3{#2}^{#3}\else#2^{#3}\fi}
	\newcommand{\id}[3][?]{#2_{#3}}
	\newcommand{\ip}[4][?]{{#2}_{#3}^{#4}}
	\newcommand{\pwR}[3][a]{\ifx!#1{\bkR[#1]{#2}}^{#3}\else\bkR[#1]{#2}^{#3}\fi}
	\newcommand{\pwB}[3][a]{\ifx!#1{\bkB[#1]{#2}}^{#3}\else\bkB[#1]{#2}^{#3}\fi}
	\newcommand{\pwS}[3][a]{\ifx!#1{\bkS[#1]{#2}}^{#3}\else\bkS[#1]{#2}^{#3}\fi}

	\newcommand{\tpT}[3][a]{ {#2}\atop \bkR[#1]{#3} }

	\newcommand{\nSmO}[2][?]{\ifx l#1\sum\limits_{#2}\else\ifx t#1{\textstyle\sum\limits_{#2}}\else\sum_{#2}\fi\fi}
	\newcommand{\nSmT}[3][?]{\ifx l#1\sum\limits_{#2}^{#3}\else\if t#1{\textstyle\sum\limits_{#2}^{#3}}\else\sum_{#2}^{#3}\fi\fi}	
	\newcommand{\nSmN}[1][?]{\ifx l#1\sum\limits\else\ifx t#1{\textstyle\sum\limits}\else\sum\fi\fi}
	\newcommand{\pSm}[2][?]{\ifx t#1 \sum_{#2}^{\prime} \else \sideset{}{^\prime}\sum_{#2} \fi}
	\newcommand{\pSmT}[3][?]{\ifx t#1 \sum_{#2}^{\prime#3} \else \sideset{}{^\prime}\sum_{#2}^{#3} \fi}	
	\newcommand{\pSmN}[1][?]{\ifx t#1 \sum^{\prime} \else \sideset{}{^\prime}\sum \fi}
	\newcommand{\dSm}[2][?]{\ifx t#1 \sum_{#2}^{\dagger} \else \sideset{}{^\dagger}\sum_{#2} \fi}
	\newcommand{\dSmT}[3][?]{\ifx t#1 \sum_{#2}^{\dagger#3} \else \sideset{}{^\dagger}\sum_{#2}^{#3} \fi}	
	\newcommand{\dSmN}[1][?]{\ifx t#1 \sum^{\dagger} \else \sideset{}{^\dagger}\sum \fi}
	\newcommand{\tpTSm}[3][?]{\nSmO[#1]{\tpT{#2}{#3}}}

		\newcommand{\Sm}{\nSmO}			\newcommand{\SmT}{\nSmT}			\newcommand{\SmN}{\nSmN}
		\newcommand{\tpSm}{\tpTSm}
		
	\newcommand{\nPd}[2][?]{\ifx l#1 \prod\limits_{#2} \else \prod_{#2} \fi}
	\newcommand{\nPdT}[3][?]{\ifx l#1 \prod\limits_{#2}^{#3} \else \prod_{#2}^{#3} \fi}

	\newcommand{\nOPs}[2][?]{\ifx l#1 \OPlus\limits_{#2} \else \OPlus_{#2} \fi}
	\newcommand{\nOPsT}[3][?]{\ifx l#1 \OPlus\limits_{#2}^{#3} \else \OPlus_{#2}^{#3} \fi}	
	\newcommand{\pOPs}[2][?]{\ifx t#1 \OPlus_{#2}^{\prime} \else \sideset{}{^\prime}\OPlus_{#2} \fi}
	\newcommand{\pOPsT}[3][?]{\ifx t#1 \OPlus_{#2}^{\prime#3} \else \sideset{}{^\prime}\OPlus_{#2}^{#3} \fi}

	\newcommand{\nIs}[2][?]{\ifx l#1 \bigcap\limits_{#2}\else\ifx b#1 \bigcap_{#2}\else{\textstyle\bigcap\limits_{#2}}\fi\fi}
	\newcommand{\nIsT}[3][?]{\ifx l#1 \bigcap\limits_{#2}^{#3}\else\ifx b#1 \bigcap_{#2}^{#3}\else{\textstyle\bigcap\limits_{#2}^{#3}}\fi\fi}	
	\newcommand{\pIs}[2][?]{\ifx t#1 \bigcap_{#2}^{\prime} \else \sideset{}{^\prime}\bigcap_{#2} \fi}
	\newcommand{\pIsT}[3][?]{\ifx t#1 \bigcap_{#2}^{\prime#3} \else \sideset{}{^\prime}\bigcap_{#2}^{#3} \fi}

	\newcommand{\nUn}[2][?]{\ifx L#1 \bigcup\limits_{#2}\else\ifx b#1 \bigcup_{#2}\else\ifx t#1{\textstyle\bigcup_{#2}}\else{\textstyle\bigcup_{#2}}\fi\fi\fi}
	\newcommand{\nUnT}[3][?]{\ifx L#1 \bigcup\limits_{#2}^{#3}\else\ifx b#1 \bigcup_{#2}^{#3}\else\ifx t#1{\textstyle\bigcup_{#2}^{#3}}\else{\textstyle\bigcup\limits_{#2}^{#3}}\fi\fi\fi}	
	\newcommand{\pUn}[2][?]{\ifx t#1 \bigcup_{#2}^{\prime} \else \sideset{}{^\prime}\bigcup_{#2} \fi}
	\newcommand{\pUnT}[3][?]{\ifx t#1 \bigcup_{#2}^{\prime#3} \else \sideset{}{^\prime}\bigcup_{#2}^{#3} \fi}	
	\newcommand{\dUn}[2][?]{\ifx L#1 \bigsqcup\limits_{#2}\else\ifx b#1 \bigsqcup_{#2}\else\ifx t#1{\textstyle\bigsqcup_{#2}}\else{\textstyle\bigsqcup\limits_{#2}}\fi\fi\fi}
	\newcommand{\dUnT}[3][?]{\ifx L#1 \bigsqcup\limits_{#2}^{#3}\else\ifx b#1 \bigsqcup_{#2}^{#3}\else\ifx t#1{\textstyle\bigsqcup_{#2}^{#3}}\else{\textstyle\bigsqcup\limits_{#2}^{#3}}\fi\fi\fi}

	\newcommand{\nLm}[2][?]{\ifx l#1 \lim\limits_{#2} \else \lim_{#2} \fi}
	\newcommand{\iLm}[2][?]{\ifx l#1 \liminf\limits_{#2} \else \liminf_{#2} \fi}
	\newcommand{\sLm}[2][?]{\ifx l#1 \limsup\limits_{#2} \else \limsup_{#2} \fi}
		
	\newcommand{\glcondEnvLineHead}[1]{ \ifx*#1 \begin{eqnarray*} \else \begin{eqnarray}  \label{#1} \fi }
	\newcommand{\glcondEnvLineTail}[1]{ \ifx*#1 \end{eqnarray*} \else \end{eqnarray} \fi }
	\newcommand{\glcondDis}[1]{\ifx d#1 \displaystyle \fi}
	\newcommand{\glcmdEqShift}{\hspace{-20pt}}
	\newcommand{\glcmdHLineCWiden}{\rule{0cm}{15pt}}	\newcommand{\glcdH}{\glcmdHLineCWiden}
	\newcommand{\lccondPar}[1]{\ifx#1p \\ \fi}

		\newcommand{\envMO}[2][*]{$\ifx d#1 \displaystyle \fi#2$}
		\newcommand{\envMT}[3][*]{$\ifx d#1 \displaystyle \fi#2=#3$}
		\newcommand{\envMTDef}[3][*]{$\ifx d#1 \displaystyle \fi#2:=#3$}
		\newcommand{\envMTPt}[4][*]{$\ifx d#1 \displaystyle \fi#3#2#4$}
			\newcommand{\envM}{\envMT}

		\newcommand{\envMTh}[4][*]{$\ifx d#1 \displaystyle \fi#2=#3=#4$}
		\newcommand{\envMThDef}[4][*]{$\ifx d#1 \displaystyle \fi#2:=#3=#4$}
		\newcommand{\envMThPt}[5][*]{$\ifx d#1 \displaystyle \fi#3#2#4#2#5$}
		\newcommand{\envMThPte}[6][*]{$\ifx d#1 \displaystyle \fi#2#3#4#5#6$}
		\newcommand{\envMF}[5][*]{$\ifx d#1 \displaystyle \fi#2=#3=#4=#5$}
		\newcommand{\envMFPt}[6][*]{$\ifx d#1 \displaystyle \fi#3#2#4#2#5#2#6$}
		
		\newcommand{\envLineT}[3][*]{ \glcondEnvLineHead{#1} & &\glcmdEqShift#2\nonumber\\&=&#3 \glcondEnvLineTail{#1} }
		
		\newcommand{\envLineTPt}[4][*]{\glcondEnvLineHead{#1} & &\glcmdEqShift#3\nonumber\\&#2&#4 \glcondEnvLineTail{#1}}
		
			\newcommand{\envLine}{\envLineT}
			
			\newcommand{\envLinePt}{\envLineTPt}
				
		\newcommand{\envHLineT}[3][*]{ \glcondEnvLineHead{#1} #2&=&#3\glcondEnvLineTail{#1} }
		
		\newcommand{\envHLineTPt}[4][*]{\glcondEnvLineHead{#1} #3&#2&#4\glcondEnvLineTail{#1}}
		
			\newcommand{\envHLine}{\envHLineT}
			
			\newcommand{\envHLinePt}{\envHLineTPt}
			
		\newcommand{\envHLineTh}[4][*]{ \glcondEnvLineHead{#1} #2&=&#3\nonumber\\&=&#4\glcondEnvLineTail{#1}}

		\newcommand{\envPLine}[2][*]{\glcondEnvLineHead{#1} #2\glcondEnvLineTail{#1}}
		
		\newcommand{\envOTLine}[4][*]{\glcondEnvLineHead{#1} #2\lnAP{=}#3\lnP{=}#4\glcondEnvLineTail{#1}}

			\newcommand{\envOTLineTh}{\envOTLine}
			
	\newcommand{\lcparaCase}{\vspace{3pt}}

	\newcommand{\envMatT}[3][a]{\autobk{(}{)}{\begin{matrix}#2\\#3\end{matrix}}}

	\newcommand{\envMatTh}[4][a]{ \autobk{(}{)}{\begin{matrix}#2\\#3\\#4\end{matrix}} }
	
	\newcommand{\envMatF}[5][a]{ \autobk{(}{)}{\begin{matrix}#2\\#3\\#4\\#5\end{matrix}} }

		\newcommand{\MatT}{\envMatT}
		
		\newcommand{\MatTh}{\envMatTh}
		
		\newcommand{\MatF}{\envMatF}

	\newcommand{\matu}[1][?]{\ifx#1?I\else I_{#1}\fi}
	\newcommand{\vlA}[2][n]{\bkAll[#1]{|}{|}{#2}}
				
	\newcommand{\nmSet}[2][n]{\vlA[#1]{#2}}	
		
	\newcommand{\sTx}[2][?]{ \ifx t#1{\tiny #2} \else \ifx s#1{\scriptsize #2} \else \ifx f#1{\footnotesize #2} \else \ifx S#1{\small #2} \else \ifx n#1{\normalsize #2} \else \ifx l#1{\large #2} \else \ifx L#1{\Large #2} \else \ifx R#1{\LARGE #2} \else \ifx h#1{\huge #2} \else \ifx H#1{\Huge #2} \else \ifx ?#1 #2 \else #2 \fi\fi\fi\fi\fi\fi\fi\fi\fi\fi\fi }
	\newcommand{\bfTx}[1]{{\bf#1}}
	
	\newcommand{\osMTx}[3][?]{\overset{#3}{#2}}
	\newcommand{\usbMTx}[3][?]{\underset{#2}{\underbrace{#3}}}

		\newcommand{\osTx}{\osMTx}
		\newcommand{\usbTx}{\usbMTx}
		
	\newcommand{\raTx}[3][?]{\raisebox{#2pt}[0pt][0pt]{\ifx d#1\displaystyle\fi#3}}
	\newcommand{\raMTx}[3][?]{\raisebox{#2pt}[0pt][0pt]{$\ifx d#1\displaystyle\fi#3$}}

	\newcommand{\roMTx}[3][?]{\rotatebox[origin=c]{#2}{$#3$}}
	
	\newcommand{\envCenter}[2][*]{\ifx*#1\begin{center}\else\begin{center}[#1]\fi #2\end{center}}
	\newcommand{\envFlushleft}[2][*]{\ifx*#1\begin{flushleft}\else\begin{flushleft}[#1]\fi #2\end{flushleft}}
	\newcommand{\envFlushright}[2][*]{\ifx*#1\begin{flushright}\else\begin{flushright}[#1]\fi #2\end{flushright}}
		
	\newcommand{\envItemIm}[2][*]{\ifx*#1\begin{itemize}\else\begin{itemize}[#1]\fi #2\end{itemize}}
	\newcommand{\envItemDp}[2][*]{\ifx*#1\begin{description}\else\begin{description}[#1]\fi #2\end{description}}
	\newcommand{\envItemEm}[2][*]{\ifx*#1\begin{enumerate}\else\begin{enfumerate}[#1]\fi #2\end{enumerate}}

	\newcommand{\envMultCol}[3][*]{\ifx1#2#3\else\begin{multicols}{#2}\ifx*#1\else\mbox{}\vspace{-#1pt}\fi#3\end{multicols}\fi}
	
\theoremstyle{plain}
\newtheorem{theorem}{THEOREM}[section]
\newtheorem{proposition}[theorem]{PROPOSITION}
\newtheorem{lemma}[theorem]{LEMMA}
\newtheorem{corollary}[theorem]{COROLLARY}
\theoremstyle{definition}

\theoremstyle{remark}
\newtheorem{remark}[theorem]{REMARK}
\theoremstyle{plain}

\theoremstyle{definition}

\theoremstyle{remark}

\theoremstyle{plain}
\newtheorem{theoremN}{THEOREM}[section]

\theoremstyle{definition}

\theoremstyle{remark}

\allowdisplaybreaks[4]
\numberwithin{equation}{section}

	\newcommand{\lccondBibitem}[3][]{ \if ?#2 \bibitem{#3} \else \bibitem[#2]{#3} \fi}
	\newcommand{\refPaper}[8][?]{
			\lccondBibitem{#1}{#2}
				#3,			
				\emph{#4}, 	
				#5\ 			
				{\bf #6}		
				(#7),			
				#8.			
		}
	\newcommand{\refPreprint}[6][?]{
			\lccondBibitem{#1}{#2}
				#3,			
				\emph{#4}, 	
				preprint; #5,	
				#6.			
		}
	\newcommand{\refPaperRep}[9][?]{
			\lccondBibitem{#1}{#2}
				#3,			
				\emph{#4}, 	
				#5\ 			
				{\bf #6}		
				(#7),			
				#8			
				; reprinted in #9	
		}

	\newcommand{\refPaperAlm}[5][?]{
			\lccondBibitem{#1}{#2}
				#3,	 		
				\emph{#4}, 	
				#5		
		}

	\newcommand{\glcondEnvLineTailPd}[1]{.\ifx*#1 \end{eqnarray*} \else \end{eqnarray} \fi  }
	\newcommand{\glcondEnvLineTailCm}[1]{,\ifx*#1 \end{eqnarray*} \else \end{eqnarray} \fi }
	\newcommand{\prcondEnvEqSpHead}[1]{ \ifx*#1 \begin{equation*}[ERROR] \else \begin{equation}  \label{#1} \fi  }
	\newcommand{\prcondEnvEqSpTail}[1]{\ifx*#1 [ERROR]\end{equation*} \else \end{equation} \fi }

	\newcommand{\envProof}[2][?]{ \par\mbox{}\vspace{-5pt}\\ \ifx?#1\emph{Proof.}\else\emph{Proof of #1.}\fi \ #2 \hfill $\Box$\\ \par}
	
		\newcommand{\envLineTPd}[3][*]{ \glcondEnvLineHead{#1} & &\glcmdEqShift#2\nonumber\\&=&#3 \glcondEnvLineTailPd{#1} }
		\newcommand{\envLineTDefPd}[3][*]{ \glcondEnvLineHead{#1} & &\glcmdEqShift#2\nonumber\\&:=&#3 \glcondEnvLineTailPd{#1} }

			\newcommand{\envLinePd}{\envLineTPd}
			\newcommand{\envLineDefPd}{\envLineTDefPd}

		\newcommand{\envLineTCmPt}[4][*]{\glcondEnvLineHead{#1} & &\glcmdEqShift#3\nonumber\\&#2&#4 \glcondEnvLineTailCm{#1}}
		\newcommand{\envLineTPdPt}[4][*]{\glcondEnvLineHead{#1} & &\glcmdEqShift#3\nonumber\\&#2&#4 \glcondEnvLineTailPd{#1}}

			\newcommand{\envLineCmPt}{\envLineTCmPt}
			\newcommand{\envLinePdPt}{\envLineTPdPt}

		\newcommand{\envLineThCmPt}[5][*]{\glcondEnvLineHead{#1}  & &\glcmdEqShift#3\nonumber\\&#2&#4\nonumber \\&#2&#5 \glcondEnvLineTailCm{#1}}

		\newcommand{\envLineFCmPte}[8][*]{\glcondEnvLineHead{#1}  & &\glcmdEqShift#2\nonumber\\&#3&#4\nonumber \\&#5&#6\nonumber \\&#7&#8 \glcondEnvLineTailCm{#1}}
		
		\newcommand{\envHLineTPd}[3][*]{ \glcondEnvLineHead{#1} #2&=&#3\glcondEnvLineTailPd{#1} }
		\newcommand{\envHLineTDefPd}[3][*]{ \glcondEnvLineHead{#1} #2&:=&#3\glcondEnvLineTailPd{#1} }
		\newcommand{\envHLineTCm}[3][*]{ \glcondEnvLineHead{#1} #2&=&#3\glcondEnvLineTailCm{#1} }
		\newcommand{\envHLineTCmDef}[3][*]{ \glcondEnvLineHead{#1} #2&:=&#3\glcondEnvLineTailCm{#1} }
		\newcommand{\envHLineTCmPt}[4][*]{\glcondEnvLineHead{#1} #3&#2&#4\glcondEnvLineTailCm{#1}}					
		\newcommand{\envHLineTPdPt}[4][*]{\glcondEnvLineHead{#1} #3&#2&#4\glcondEnvLineTailPd{#1}}

			\newcommand{\envHLinePd}{\envHLineTPd}
			\newcommand{\envHLineDefPd}{\envHLineTDefPd}
			\newcommand{\envHLineCm}{\envHLineTCm}
			\newcommand{\envHLineCmDef}{\envHLineTCmDef}
			\newcommand{\envHLineCmPt}{\envHLineTCmPt}
			\newcommand{\envHLinePdPt}{\envHLineTPdPt}

		\newcommand{\envHLineThPd}[4][*]{ \glcondEnvLineHead{#1} #2&=&#3\nonumber\\&=&#4\glcondEnvLineTailPd{#1}	 }
		
		\newcommand{\envHLineThCm}[4][*]{ \glcondEnvLineHead{#1} #2&=&#3\nonumber\\&=&#4\glcondEnvLineTailCm{#1}}
		
		\newcommand{\envHLineThCmPt}[5][*]{\glcondEnvLineHead{#1} #3&#2&#4\nonumber\\&#2&#5  \glcondEnvLineTailCm{#1}}

		\newcommand{\envHLineFCmPt}[6][*]{\glcondEnvLineHead{#1} #3&#2&#4\nonumber\\&#2&#5\nonumber \\&#2&#6\glcondEnvLineTailCm{#1}}

		\newcommand{\envHLineCFCmDef}[5][*]{\glcondEnvLineHead{#1} #2&:=&#3,\nonumber\\\glcdH#4&:=&#5\glcondEnvLineTailCm{#1}}
		
		\newcommand{\envHLineCFCmNme}[5][*]{\begin{eqnarray} #2&=&#3,\\\glcdH#4&=&#5 \glcondEnvLineTailCm{?} }
		\newcommand{\envHLineCFNmePd}[5][*]{\begin{eqnarray} #2&=&#3,\\\glcdH#4&=&#5 \glcondEnvLineTailPd{?} }
		\newcommand{\envHLineCFCmDefNme}[5][*]{\begin{eqnarray} #2&:=&#3,\\\glcdH#4&:=&#5 \glcondEnvLineTailCm{?} }
		\newcommand{\envHLineCFDefNmePd}[5][*]{\begin{eqnarray} #2&:=&#3,\\\glcdH#4&:=&#5 \glcondEnvLineTailPd{?} }

		\newcommand{\envHLineCFNmePdPt}[6][*]{\begin{eqnarray}#3&#2&#4,\\\glcdH#5&#2&#6\glcondEnvLineTailPd{?}}
		\newcommand{\envHLineCFCmNmePt}[6][*]{\begin{eqnarray}#3&#2&#4,\\\glcdH#5&#2&#6\glcondEnvLineTailCm{?}}
		\newcommand{\envHLineCFNmePdPte}[7][*]{\begin{eqnarray}#2&#3&#4,\\\glcdH#5&#6&#7\glcondEnvLineTailPd{?}}
		\newcommand{\envHLineCFCmNmePte}[7][*]{\begin{eqnarray}#2&#3&#4,\\\glcdH#5&#6&#7\glcondEnvLineTailCm{?}}

		\newcommand{\envHLineCFLqqPd}[5][*]{\envPLinePd[#1]{#2\lnP{=}#3,\qquad#4\lnP{=}#5}}

		\newcommand{\envHLineCFLaaCm}[5][*]{\envPLineCm[#1]{#2\lnP{=}#3\qquad\text{and}\qquad#4\lnP{=}#5}}

		\newcommand{\envHLineCSPd}[7][*]{\glcondEnvLineHead{#1} #2&=&#3,\nonumber\\\glcdH#4&=&#5,\\\glcdH#6&=&#7\nonumber\glcondEnvLineTailPd{#1}}
		
		\newcommand{\envHLineCSCm}[7][*]{\glcondEnvLineHead{#1} #2&=&#3,\nonumber\\\glcdH#4&=&#5,\\\glcdH#6&=&#7\nonumber\glcondEnvLineTailCm{#1}}
			
		\newcommand{\envHLineCSNmePd}[7][*]{\begin{eqnarray} #2&=&#3,\\\glcdH#4&=&#5,\\\glcdH#6&=&#7\glcondEnvLineTailPd{?}}
		\newcommand{\envHLineCSDefNmePd}[7][*]{\begin{eqnarray} #2&:=&#3,\\\glcdH#4&:=&#5,\\\glcdH#6&:=&#7\glcondEnvLineTailPd{?}}
		\newcommand{\envHLineCSCmNme}[7][*]{\begin{eqnarray} #2&=&#3,\\\glcdH#4&=&#5,\\\glcdH#6&=&#7\glcondEnvLineTailCm{?}}
		\newcommand{\envHLineCSCmDefNme}[7][*]{\begin{eqnarray} #2&:=&#3,\\\glcdH#4&:=&#5,\\\glcdH#6&:=&#7\glcondEnvLineTailCm{?}}

		\newcommand{\envHLineCSNmePdPt}[8][*]{\begin{eqnarray}#3&#2&#4,\\\glcdH#5&#2&#6,\\\glcdH#7&#2&#8\glcondEnvLineTailPd{?}}
		\newcommand{\envHLineCSCmNmePt}[8][*]{\begin{eqnarray}#3&#2&#4,\\\glcdH#5&#2&#6,\\\glcdH#7&#2&#8\glcondEnvLineTailCm{?}}
		\newcommand{\envHLineCSNmePdPte}[9][*]{\begin{eqnarray}#2&#3&#4,\\\glcdH#5&#6&#7,\\\glcdH#8&#2&#9\glcondEnvLineTailPd{?}}
		\newcommand{\envHLineCSCmNmePte}[9][*]{\begin{eqnarray}#2&#3&#4,\\\glcdH#5&#6&#7,\\\glcdH#8&#2&#9\glcondEnvLineTailCm{?}}
		
		\newcommand{\envHLineCECm}[9][*]{\glcondEnvLineHead{#1} #2&=&#3,\nonumber\\\glcdH#4&=&#5,\nonumber\\\glcdH#6&=&#7,\\\glcdH#8&=&#9\nonumber\glcondEnvLineTailCm{#1}}
		
		\newcommand{\envHLineCENmePd}[9][*]{\begin{eqnarray} #2&=&#3,\\\glcdH#4&=&#5,\\\glcdH#6&=&#7,\\\glcdH#8&=&#9\glcondEnvLineTailPd{?}}
		\newcommand{\envHLineCEDefNmePd}[9][*]{\begin{eqnarray} #2&:=&#3,\\\glcdH#4&:=&#5,\\\glcdH#6&:=&#7,\\\glcdH#8&:=&#9\glcondEnvLineTailPd{?}}
		\newcommand{\envHLineCECmNme}[9][*]{\begin{eqnarray} #2&=&#3,\\\glcdH#4&=&#5,\\\glcdH#6&=&#7,\\\glcdH#8&=&#9\glcondEnvLineTailCm{?}}
		\newcommand{\envHLineCECmDefNme}[9][*]{\begin{eqnarray} #2&:=&#3,\\\glcdH#4&:=&#5,\\\glcdH#6&:=&#7,\\\glcdH#8&:=&#9\glcondEnvLineTailCm{?}}

			\newcommand{\pccondPaForPar}[1]{\ifx#1p \\\glcdH \fi}
			\newcommand{\pccondPaForNonnum}[1]{\ifx#1* \nonumber \fi}

		\newcommand{\envPLinePd}[2][*]{\glcondEnvLineHead{#1} #2\glcondEnvLineTailPd{#1}}
		\newcommand{\envPLineCm}[2][*]{\glcondEnvLineHead{#1} #2\glcondEnvLineTailCm{#1}}
	
		\newcommand{\envOTLineCm}[4][*]{\glcondEnvLineHead{#1} #2\lnAP{=}#3\lnP{=}#4,\glcondEnvLineTail{#1}}
		\newcommand{\envOTLinePd}[4][*]{\glcondEnvLineHead{#1} #2\lnAP{=}#3\lnP{=}#4.\glcondEnvLineTail{#1}}
		
		\newcommand{\envOTLineDefPd}[4][*]{\glcondEnvLineHead{#1} #2\lnAP{:=}#3\lnP{=}#4.\glcondEnvLineTail{#1}}

			\newcommand{\envOTLineThCm}{\envOTLineCm}

		\newcommand{\envOFLineCm}[5][*]{\glcondEnvLineHead{#1} #2\lnAP{=}#3\lnP{=}#4\lnP{=}#5,\glcondEnvLineTail{#1}}

			\newcommand{\envOFLineFCm}{\envOFLineCm}

		\newcommand{\envMOCm}[2][*]{$\ifx d#1 \displaystyle \fi#2$,}
		\newcommand{\envMOPd}[2][*]{$\ifx d#1 \displaystyle \fi#2$.}
		
		\newcommand{\envMTCm}[3][*]{$\ifx d#1 \displaystyle \fi#2=#3$,}
		\newcommand{\envMTPd}[3][*]{$\ifx d#1 \displaystyle \fi#2=#3$.}
		\newcommand{\envMTCmDef}[3][*]{$\ifx d#1 \displaystyle \fi#2:=#3$,}
		\newcommand{\envMTDefPd}[3][*]{$\ifx d#1 \displaystyle \fi#2:=#3$.}
		\newcommand{\envMTCmPt}[4][*]{$\ifx d#1 \displaystyle \fi#3#2#4$,}
		\newcommand{\envMTPdPt}[4][*]{$\ifx d#1 \displaystyle \fi#3#2#4$.}
			\newcommand{\envMCm}{\envMTCm}
			\newcommand{\envMPd}{\envMTPd}
			\newcommand{\envMCmDef}{\envMTCmDef}

		\newcommand{\envMThCm}[4][*]{$\ifx d#1 \displaystyle \fi#2=#3=#4$,}
		\newcommand{\envMThPd}[4][*]{$\ifx d#1 \displaystyle \fi#2=#3=#4$.}
		\newcommand{\envMThCmPt}[5][*]{$\ifx d#1 \displaystyle \fi#3#2#4#2#5$,}
		\newcommand{\envMThPdPt}[5][*]{$\ifx d#1 \displaystyle \fi#3#2#4#2#5$.}
		\newcommand{\envMFCm}[5][*]{$\ifx d#1 \displaystyle \fi#2=#3=#4=#5$,}
		\newcommand{\envMFPd}[5][*]{$\ifx d#1 \displaystyle \fi#2=#3=#4=#5$.}
		\newcommand{\envMFCmPt}[6][*]{$\ifx d#1 \displaystyle \fi#3#2#4#2#5#2#6$,}
		\newcommand{\envMFPdPt}[6][*]{$\ifx d#1 \displaystyle \fi#3#2#4#2#5#2#6$.}

		\newcommand{\prcondHLCPNm}{\hspace{-1pt}}
		\newcommand{\envHLineCFCmNm}[5][*]{ \begin{equation}\begin{split} \ifx*#1 \text{[ERROR;need label name]} \else \label{#1} \fi #2&\prcondHLCPNm\lnP{=}\prcondHLCPNm#3,\\#4&\prcondHLCPNm\lnP{=}\prcondHLCPNm#5, \end{split}\end{equation} }
		\newcommand{\envHLineCFNm}[5][*]{ \begin{equation}\begin{split} \ifx*#1 \text{[ERROR;need label name]} \else \label{#1} \fi #2&\prcondHLCPNm\lnP{=}\prcondHLCPNm#3\\#4&\prcondHLCPNm\lnP{=}\prcondHLCPNm#5, \end{split}\end{equation} }
		\newcommand{\envHLineCFNmPd}[5][*]{ \begin{equation}\begin{split} \ifx*#1 \text{[ERROR;need label name]} \else \label{#1} \fi #2&\prcondHLCPNm\lnP{=}\prcondHLCPNm#3,\\#4&\prcondHLCPNm\lnP{=}\prcondHLCPNm#5. \end{split}\end{equation} }
		\newcommand{\envHLineCFCmDefNm}[5][*]{ \begin{equation}\begin{split} \ifx*#1 \text{[ERROR;need label name]} \else \label{#1} \fi #2&\prcondHLCPNm\lnP{:=}\prcondHLCPNm#3,\\#4&\prcondHLCPNm\lnP{:=}\prcondHLCPNm#5, \end{split}\end{equation} }
		\newcommand{\envHLineCFDefNm}[5][*]{ \begin{equation}\begin{split} \ifx*#1 \text{[ERROR;need label name]} \else \label{#1} \fi #2&\prcondHLCPNm\lnP{:=}\prcondHLCPNm#3\\#4&\prcondHLCPNm\lnP{:=}\prcondHLCPNm#5, \end{split}\end{equation} }
		\newcommand{\envHLineCFDefNmPd}[5][*]{ \begin{equation}\begin{split} \ifx*#1 \text{[ERROR;need label name]} \else \label{#1} \fi #2&\prcondHLCPNm\lnP{:=}\prcondHLCPNm#3,\\#4&\prcondHLCPNm\lnP{:=}\prcondHLCPNm#5. \end{split}\end{equation} }
		\newcommand{\envHLineCSCmNm}[7][*]{ \begin{equation}\begin{split} \ifx*#1 \text{[ERROR;need label name]} \else \label{#1} \fi #2&\prcondHLCPNm\lnP{=}\prcondHLCPNm#3,\\#4&\prcondHLCPNm\lnP{=}\prcondHLCPNm#5,\\#6&\prcondHLCPNm\lnP{=}\prcondHLCPNm#7 \end{split}\end{equation} }
		\newcommand{\envHLineCSNm}[7][*]{ \begin{equation}\begin{split} \ifx*#1 \text{[ERROR;need label name]} \else \label{#1} \fi #2&\prcondHLCPNm\lnP{=}\prcondHLCPNm#3\\#4&\prcondHLCPNm\lnP{=}\prcondHLCPNm#5\\#6&\prcondHLCPNm\lnP{=}\prcondHLCPNm#7 \end{split}\end{equation} }
		\newcommand{\envHLineCSNmPd}[7][*]{ \begin{equation}\begin{split} \ifx*#1 \text{[ERROR;need label name]} \else \label{#1} \fi #2&\prcondHLCPNm\lnP{=}\prcondHLCPNm#3,\\#4&\prcondHLCPNm\lnP{=}\prcondHLCPNm#5,\\#6&\prcondHLCPNm\lnP{=}\prcondHLCPNm#7. \end{split}\end{equation} }
		\newcommand{\envHLineCSCmDefNm}[7][*]{ \begin{equation}\begin{split} \ifx*#1 \text{[ERROR;need label name]} \else \label{#1} \fi #2&\prcondHLCPNm\lnP{:=}\prcondHLCPNm#3,\\#4&\prcondHLCPNm\lnP{:=}\prcondHLCPNm#5,\\#6&\prcondHLCPNm\lnP{:=}\prcondHLCPNm#7 \end{split}\end{equation} }
		\newcommand{\envHLineCSDefNm}[7][*]{ \begin{equation}\begin{split} \ifx*#1 \text{[ERROR;need label name]} \else \label{#1} \fi #2&\prcondHLCPNm\lnP{:=}\prcondHLCPNm#3\\#4&\prcondHLCPNm\lnP{:=}\prcondHLCPNm#5\\#6&\prcondHLCPNm\lnP{:=}\prcondHLCPNm#7 \end{split}\end{equation} }
		\newcommand{\envHLineCSDefNmPd}[7][*]{ \begin{equation}\begin{split} \ifx*#1 \text{[ERROR;need label name]} \else \label{#1} \fi #2&\prcondHLCPNm\lnP{:=}\prcondHLCPNm#3,\\#4&\prcondHLCPNm\lnP{:=}\prcondHLCPNm#5,\\#6&\prcondHLCPNm\lnP{:=}\prcondHLCPNm#7. \end{split}\end{equation} }
		\newcommand{\envHLineCECmNm}[9][*]{ \begin{equation}\begin{split} \ifx*#1 \text{[ERROR;need label name]} \else \label{#1} \fi #2&\prcondHLCPNm\lnP{=}\prcondHLCPNm#3,\\#4&\prcondHLCPNm\lnP{=}\prcondHLCPNm#5,\\#6&\prcondHLCPNm\lnP{=}\prcondHLCPNm#7,\\#8&\prcondHLCPNm\lnP{=}\prcondHLCPNm#9,  \end{split}\end{equation} }
		\newcommand{\envHLineCENm}[9][*]{ \begin{equation}\begin{split} \ifx*#1 \text{[ERROR;need label name]} \else \label{#1} \fi #2&\prcondHLCPNm\lnP{=}\prcondHLCPNm#3\\#4&\prcondHLCPNm\lnP{=}\prcondHLCPNm#5\\#6&\prcondHLCPNm\lnP{=}\prcondHLCPNm#7\\#8&\prcondHLCPNm\lnP{=}\prcondHLCPNm#9  \end{split}\end{equation} }
		\newcommand{\envHLineCENmPd}[9][*]{ \begin{equation}\begin{split} \ifx*#1 \text{[ERROR;need label name]} \else \label{#1} \fi #2&\prcondHLCPNm\lnP{=}\prcondHLCPNm#3,\\#4&\prcondHLCPNm\lnP{=}\prcondHLCPNm#5,\\#6&\prcondHLCPNm\lnP{=}\prcondHLCPNm#7,\\#8&\prcondHLCPNm\lnP{=}\prcondHLCPNm#9.  \end{split}\end{equation} }
		\newcommand{\envHLineCECmDefNm}[9][*]{ \begin{equation}\begin{split} \ifx*#1 \text{[ERROR;need label name]} \else \label{#1} \fi #2&\prcondHLCPNm\lnP{:=}\prcondHLCPNm#3,\\#4&\prcondHLCPNm\lnP{:=}\prcondHLCPNm#5,\\#6&\prcondHLCPNm\lnP{:=}\prcondHLCPNm#7,\\#8&\prcondHLCPNm\lnP{:=}\prcondHLCPNm#9,  \end{split}\end{equation} }
		\newcommand{\envHLineCEDefNm}[9][*]{ \begin{equation}\begin{split} \ifx*#1 \text{[ERROR;need label name]} \else \label{#1} \fi #2&\prcondHLCPNm\lnP{:=}\prcondHLCPNm#3\\#4&\prcondHLCPNm\lnP{:=}\prcondHLCPNm#5\\#6&\prcondHLCPNm\lnP{:=}\prcondHLCPNm#7\\#8&\prcondHLCPNm\lnP{:=}\prcondHLCPNm#9  \end{split}\end{equation} }
		\newcommand{\envHLineCEDefNmPd}[9][*]{ \begin{equation}\begin{split} \ifx*#1 \text{[ERROR;need label name]} \else \label{#1} \fi #2&\prcondHLCPNm\lnP{:=}\prcondHLCPNm#3,\\#4&\prcondHLCPNm\lnP{:=}\prcondHLCPNm#5,\\#6&\prcondHLCPNm\lnP{:=}\prcondHLCPNm#7,\\#8&\prcondHLCPNm\lnP{:=}\prcondHLCPNm#9.  \end{split}\end{equation} }
	\newcommand{\envMLineTPd}[3][*]{ \ifx*#1 \begin{multline*} #2\lnP{=}#3.\end{multline*} \else \begin{multline} \label{#1} #2\lnP{=}#3.\end{multline} \fi }
	\newcommand{\envMLineTCm}[3][*]{ \ifx*#1 \begin{multline*} #2\lnP{=}#3,\end{multline*} \else \begin{multline} \label{#1} #2\lnP{=}#3,\end{multline} \fi }
	\newcommand{\envMLineTDefPd}[3][*]{ \ifx*#1 \begin{multline*} #2\lnP{:=}#3.\end{multline*} \else \begin{multline} \label{#1} #2\lnP{:=}#3.\end{multline} \fi }
	\newcommand{\envMLineTCmDef}[3][*]{ \ifx*#1 \begin{multline*} #2\lnP{:=}#3,\end{multline*} \else \begin{multline} \label{#1} #2\lnP{:=}#3,\end{multline} \fi }

	\newcommand{\envCaseFCm}[5][?]{\begin{cases} \glcondDis{#1}#2,\lcparaCase\\\glcondDis{#1}#3,\lcparaCase\\\glcondDis{#1}#4,\lcparaCase\\\glcondDis{#1}#5,\end{cases}}
	\newcommand{\envCaseFPd}[5][?]{\begin{cases} \glcondDis{#1}#2,\lcparaCase\\\glcondDis{#1}#3,\lcparaCase\\\glcondDis{#1}#4,\lcparaCase\\\glcondDis{#1}#5.\end{cases}}
	\newcommand{\envCaseFiCm}[6][?]{\begin{cases} \glcondDis{#1}#2,\lcparaCase\\\glcondDis{#1}#3,\lcparaCase\\\glcondDis{#1}#4,\lcparaCase\\\glcondDis{#1}#5,\lcparaCase\\\glcondDis{#1}#6,\end{cases}}
	\newcommand{\envCaseFiPd}[6][?]{\begin{cases} \glcondDis{#1}#2,\lcparaCase\\\glcondDis{#1}#3,\lcparaCase\\\glcondDis{#1}#4,\lcparaCase\\\glcondDis{#1}#5,\lcparaCase\\\glcondDis{#1}#6.\end{cases}}

	\DeclareFontEncoding{OT2}{}{}
	\DeclareFontSubstitution{OT2}{cmr}{m}{n}
	\DeclareFontFamily{OT2}{cmr}{\hyphenchar\font45}
	\DeclareFontShape{OT2}{cmr}{m}{n}{<5><6><7><8><9>gen*wncyr <10><10.95><12><14.4><17.28><20.74><24.88>wncyr10}{}
	\DeclareFontShape{OT2}{cmr}{b}{n}{<5><6><7><8><9>gen*wncyb<10><10.95><12><14.4><17.28><20.74><24.88>wncyb10}{}
	\DeclareMathAlphabet{\mathcyr}{OT2}{cmr}{m}{n}
	\DeclareMathAlphabet{\mathcyb}{OT2}{cmr}{b}{n}
	\SetMathAlphabet{\mathcyr}{bold}{OT2}{cmr}{b}{n}

	\newcommand{\shs}{\mathcyr{sh}}

		\newcommand{\sh}{\shs}
	\newcommand{\ztRLettHelp}{\bullet}
	\newcommand{\ztN}[1][?]{\zeta}				\newcommand{\ztO}[2][n]{\Fc[#1]{\zeta}{#2}}							
		\newcommand{\zt}{\ztO}
	\newcommand{\ztHN}[1][?]{\zeta_*}											
		
	\newcommand{\ztRN}[1][?]{\zeta_\ztRLettHelp}						
					
	\newcommand{\ztSN}[1][?]{\zeta_\sh}			\newcommand{\ztSO}[2][n]{\idFc[#1]{\zeta}{\sh}{#2}}					
		\newcommand{\ztS}{\ztSO}
	\newcommand{\ztsLettHelp}{\star}
	\newcommand{\ztsN}[1][?]{\zeta^\ztsLettHelp}											
		
	\newcommand{\ztsHN}[1][?]{\zeta_*^\ztsLettHelp}											
		
	\newcommand{\ztsRN}[1][?]{\zeta_\ztRLettHelp^\ztsLettHelp}					
					
	\newcommand{\ztsSN}[1][?]{\zeta_\sh^\ztsLettHelp}

	\newcommand{\pVMpT}[3][?]{{#2}_{[#3]}}		\newcommand{\pVMp}{\pVMpT}		\newcommand{\VMp}{\pVMp}

	\newcommand{\pmdE}[1][?]{\SpO{E_0}}
		\newcommand{\mdE}{\pmdE}
	
	\newcommand{\LetHelpFcCh}{\chi}
	\newcommand{\fcMChSN}[1][?]{\LetHelpFcCh_{\sh}}						
		
	\newcommand{\pveEO}[2][?]{\bfTx{e}_{#2}}		\newcommand{\pveE}{\pveEO}		\newcommand{\veE}{\pveE}
	
	\newcommand{\prgM}[2]{M_{#1}(#2)}
	\newcommand{\prgFS}[3][?]{#2\bkSd[n]{#3}}	
		\newcommand{\rgM}{\prgM}		
		\newcommand{\rgFS}{\prgFS}
		
	\newcommand{\pdpt}[2][n]{\Fc[#1]{\mathrm{d}}{#2}}
	\newcommand{\pwgt}[2][n]{\Fc[#1]{\mathrm{w}}{#2}}
		\newcommand{\dpt}{\pdpt}		
		\newcommand{\wgt}{\pwgt}
	\newcommand{\pbfADM}{{\bf adm}}
	
		\newcommand{\bfADM}{\pbfADM}

	\newcommand{\pmpSN}{\Sig}	\newcommand{\pmpSO}[2][?]{\pmpSN_{#2}}
		\newcommand{\mpS}{\pmpSO}		
	
	\newcommand{\pgpC}[2][?]{\mathfrak{C}_{#2}}			\newcommand{\gpC}{\pgpC}
	\newcommand{\pgpCs}[2][?]{\mathfrak{C}_{#2}^0}		\newcommand{\gpCs}{\pgpCs}
	\newcommand{\pgpST}[1][?]{\mathfrak{S}'_{3}}		\newcommand{\gpST}{\pgpST}
	\newcommand{\pgeTO}[2][?]{t_{#2}}						
	\newcommand{\pgeS}[1][?]{\ifx?#1s\else s_{#1}\fi}		\newcommand{\pgeSt}[1][?]{\ifx?#1s^{t}\else s_{#1}^{t}\fi}		\newcommand{\pgeSi}[1][?]{\ifx?#1s^{-1}\else s_{#1}^{-1}\fi}		
	\newcommand{\pgeR}[1][?]{\ifx?#1r\else r_{#1}\fi}		\newcommand{\pgeRt}[1][?]{\ifx?#1r^{t}\else r_{#1}^{t}\fi}		\newcommand{\pgeRi}[1][?]{\ifx?#1r^{-1}\else r_{#1}^{-1}\fi}	
	\newcommand{\pgeQ}[1][?]{\ifx?#1q\else q_{#1}\fi}		\newcommand{\pgeQt}[1][?]{\ifx?#1q^{t}\else q_{#1}^{t}\fi}		\newcommand{\pgeQi}[1][?]{\ifx?#1q^{-1}\else q_{#1}^{-1}\fi}
	\newcommand{\pgeP}[1][?]{\ifx?#1p\else p_{#1}\fi}		\newcommand{\pgePt}[1][?]{\ifx?#1p^{t}\else p_{#1}^{t}\fi}		\newcommand{\pgePi}[1][?]{\ifx?#1p^{-1}\else p_{#1}^{-1}\fi}	
		\newcommand{\geT}{\pgeTO}		
		\newcommand{\geS}{\pgeS}				
		\newcommand{\geR}{\pgeR}				
		\newcommand{\geQ}{\pgeQ}					
		\newcommand{\geP}{\pgeP}		\newcommand{\gePt}{\pgePt}		\newcommand{\gePi}{\pgePi}	
	\newcommand{\pgeJ}[1][?]{\ifx?#1j\else j_{#1}\fi}
		\newcommand{\geJ}{\pgeJ}
	\newcommand{\setV}{V}
	\newcommand{\pgePh}[1][?]{\ifx?#1\Phi\else\Phi_{#1}\fi}		\newcommand{\pgePs}[1][?]{\ifx?#1\Psi\else\Psi_{#1}\fi}		\newcommand{\pgeOm}[1][?]{\ifx?#1\Ome\else\Ome_{#1}\fi}
		\newcommand{\gePh}{\pgePh}		\newcommand{\gePs}{\pgePs}		\newcommand{\geOm}{\pgeOm}

	\newcommand{\pgeG}[1][?]{\ifx?#1g\else g_{#1}\fi}		\newcommand{\pgeGt}[1][?]{\ifx?#1g^{t}\else g_{#1}^{t}\fi}		
						
	\renewcommand{\gpSym}[2][?]{\mathfrak{S}_{#2}}
	\renewcommand{\gpAlt}[2][?]{\mathfrak{A}_{#2}}

	\newcommand{\pggfcQZSj}[2][n]{\idFc[#1]{\mLt{\ggfcQZlett}}{\sh}{#2}}		\newcommand{\ggfcQZSj}{\pggfcQZSj}
			\newcommand{\ggfcSZlett}{\mathcal{{S}}}					

	\newcommand{\ggfcSZ}[2][n]{\Fc[#1]{\ggfcSZlett}{#2}}								\newcommand{\ggfcSZS}[2][n]{\idFc[#1]{\ggfcSZlett}{\sh}{#2}}

			\newcommand{\ggfcDZlett}{\mathcal{{D}}}					

	\newcommand{\ggfcDZ}[2][n]{\Fc[#1]{\ggfcDZlett}{#2}}			\newcommand{\ggfcDZS}[2][n]{\idFc[#1]{\ggfcDZlett}{\sh}{#2}}

			\newcommand{\ggfcTZlett}{\mathcal{{T}}}					

	\newcommand{\ggfcTZ}[2][n]{\Fc[#1]{\ggfcTZlett}{#2}}			\newcommand{\ggfcTZS}[2][n]{\idFc[#1]{\ggfcTZlett}{\sh}{#2}}

	\newcommand{\gfcQZlett}{\mathcal{Q}}		\newcommand{\ggfcQZlett}{\mathcal{{Q}}}					
	\newcommand{\gfcQZ}[3][n]{\idFc[#1]{\gfcQZlett}{#2}{#3}}					
	\newcommand{\ggfcQZ}[2][n]{\Fc[#1]{\ggfcQZlett}{#2}}			\newcommand{\ggfcQZS}[2][n]{\idFc[#1]{\ggfcQZlett}{\sh}{#2}}

	\newcommand{\pggfcSdZ}[2][n]{\Fc[#1]{{\mathfrak{d}}}{#2}}					\newcommand{\ggfcSdZ}{\pggfcSdZ}	
	\newcommand{\pggfcStZ}[2][n]{\Fc[#1]{{\mathfrak{t}}}{#2}}						\newcommand{\ggfcStZ}{\pggfcStZ}	
	\newcommand{\pggfcSqZ}[2][n]{\Fc[#1]{{\mathfrak{q}}}{#2}}					\newcommand{\ggfcSqZ}{\pggfcSqZ}		
	\newcommand{\pgfcSqZ}[3][n]{\idFc[#1]{{\mathfrak{q}}}{#2}{#3}}					\newcommand{\gfcSqZ}{\pgfcSqZ}	
			\newcommand{\ggfcZlett}{\mathcal{{Z}}}					
	\newcommand{\ggfcZ}[3][n]{\idFc[#1]{\ggfcZlett}{#2}{#3}}						
	\newcommand{\ggfcZS}[3][n]{\idFc[#1]{\ggfcZlett}{\sh,#2}{#3}}					
	\newcommand{\popSPlett}{\sharp}
	\newcommand{\popSP}[2][?]{\bkR[#1]{#2}^\popSPlett}				\newcommand{\popSPr}[2][n]{(#2)^\popSPlett}
	\newcommand{\popSPT}[3][n]{\pwFc[#1]{#2}{\popSPlett}{#3}}			
	\newcommand{\lfgSHiN}{sh}		\newcommand{\lfgSHiT}[3][?]{\ifx?#1\lfgSHiN_{#2}^{(#3)}\else\pLt{\lfgSHiN}_{#2}^{(#3)}\fi}	
								\newcommand{\lfgSHiO}[2][?]{\lfgSHiN_{#2}}		\newcommand{\lfgSHi}{\lfgSHiO}
	\newcommand{\popR}[3][n]{\idFc[#1]{\mathbf{r}}{#2}{#3}}
		\newcommand{\opR}{\popR}	

\allowdisplaybreaks[4]
	\title{\mainTitle}
	\author{\authorName
			\thanks{\organizationNameFst, \placeAddressFst}
			\mbox{}
			\thanks{\organizationNameScd, \departmentNameScd} 
		}
	\date{}

\begin{document}
\maketitle
\renewcommand{\thefootnote}{\fnsymbol{footnote}}
\footnote[0]{e-mail : \emailAddressFst}
\footnote[0]{MSC-class: \MSCname}
\footnote[0]{Key words: \keyWord}
\renewcommand{\thefootnote}{\arabic{footnote}}\setcounter{footnote}{0}
\vPack[30]

\begin{abstract}
In the present paper,
	we prove an identity for the generating function of the quadruple zeta values with the action of the matrix ring $\rgM{4}{\setZ}$.
Taking homogeneous parts on both sides of the identity and substituting appropriate values for the variables, 
	we obtain the sum formula for quadruple zeta values.
We also obtain its weighted analogues,
	which include the formulas for this case proved by 
	Guo and Xie (2009, J. Number Theory 129, 2747--2765) 
	and by Ong, Eie, and Liaw (2013, Int. J. Number Theory 9, 1185--1198).
\end{abstract}

\section{Introduction and statement of results} \label{sectOne}
A multiple zeta value (MZV) is defined by the convergent series
	\envOTLineCm
	{
		\zt{\bfL_r}
	}
	{
		\zt{l_1,l_2,\ldots,l_r}
	}
	{
		\Sm{m_1>m_2>\cdots>m_r>0} \frc{ \pw{m_1}{l_1}\pw{m_2}{l_2}\cdots\pw{m_r}{l_r} }{1}
	}
	where $\bfL_r=(l_1,l_2,\ldots,l_r)$ is an admissible index set,
	that is, 
	it is a sequence of positive integers with $l_1\geq2$.
The condition $l_1\geq2$ ensures convergence.
The integers $\dpt{\bfL_r}=r$ and $\wgt{\bfL_r}=l_1+\cdots+l_r$ are called the depth and weight, 
	respectively.
The single, double, triple, and quadruple zeta values (SZVs, DZVs, TZVs, and QZVs, respectively)
	are the MZVs of depth $1$, $2$, $3$, and $4$, respectively.

It is known that, among the MZVs, there are many linear relations over $\setZ$.
One notable example is the sum formula,
	which was conjectured by Moen and Markett (see \cite{Hoffman92} and \cite{Markett94}, respectively):
	\envHLineCm[1_PL_EqSumFml]
	{
		\tpSm{ \bfL_r : \bfADM }{ \wgt{\bfL_r}=l } \zt{l_1,\ldots, l_r}
	}
	{
		\zt{l}
	}
	where the summation ranges over all admissible index sets of depth $r$ and weight $l$.
Formula \refEq{1_PL_EqSumFml} was proved for DZVs by Euler \cite{Euler1776}, for TZVs by Hoffman and Moen \cite{HM96},
	and for general MZVs by Granville \cite{Granville97}.
Zagier also proved \refEq{1_PL_EqSumFml} independently in an unpublished manuscript (see \cite{Granville97}).

Formula \refEq{1_PL_EqSumFml} has been generalized and extended in various directions 
	\cite{ELO09,GX09,Ho12Ax,HO03,Machide12,Machide13a,Nakamura09,Ohno99,OZu08,OEL13,SC12}.
Recently,
	generalizations for DZVs and TZVs, from the point of view of the generating functions, 
	were given in \cite{GKZ06} and \cite{Machide13b},
	respectively (see \ref{sectAppA} for details).
In the present paper,
	we give such a generalization of \refEq{1_PL_EqSumFml} for QZVs.

We will begin by introducing the notation and terminology that will be used to state our results.
We define the generating function $\ggfcQZ[]{}$ of the QZVs as the formal power series, 
	\envHLineDefPd[1_PL_DefGGfcQZV]
	{
		\ggfcQZ{x_1,x_2,x_3,x_4}
	}
	{
		\Sm{ \bfL_4 : \bfADM } \zt{\bfL_4} x_1^{l_1-1} x_2^{l_2-1} x_3^{l_3-1} x_4^{l_4-1}
	}	
Replacing $\zt{\bfL_4}$ with $\zt{ \wgt{\bfL_4} } = \zt{l_1+l_2+l_3+l_4}$,
	we also define 
	\envHLineCmDef[1_PL_DefFcSqZ]
	{
		\ggfcSqZ{x_1,x_2,x_3,x_4}
	}
	{
		\tpSm{\bfL_4 }{ \wgt{\bfL_4}>4 } \zt{ \wgt{\bfL_4} } x_1^{l_1-1} x_2^{l_2-1} x_3^{l_3-1} x_4^{l_4-1}
	}
	where the summation rules are different:
	the rule of \refEq{1_PL_DefGGfcQZV} ranges over all admissible index sets whose depth is $4$,
	but that of \refEq{1_PL_DefFcSqZ} ranges over all sequences of positive integers whose depth is $4$ and whose weight is greater than $4$.

We denote by $\gpSym{r}$ and $\gpu=\gpu[r]$ the symmetric group of degree $r$ and its identity element,
	respectively,
	and 
	also denote by $\gpGL{r}{\setZ}$ and $\matI=\matI[r]$ the general linear group of degree $r$ over $\setZ$ and its identity element,
	respectively.
We identify a permutation $\sig$ in $\gpSym{r}$ with the matrix $(\delta_{i\sig(j)})_{1\leq i,j\leq r}$ in $\gpGL{r}{\setZ}$,
	where $\delta_{ij}$ is the Kronecker delta function.
For example,
	\envPLinePd
	{
		\gpu		\lnP{=}	\matI
	,\qquad
		(1234)		\lnP{=}	\MatF{0&0&0&1}{1&0&0&0}{0&1&0&0}{0&0&1&0}		
	,\quad\text{and}\quad
		(13)(24)		\lnP{=}	\MatF{0&0&1&0}{0&0&0&1}{1&0&0&0}{0&1&0&0}		
	}
Let $\mopTP{M}$ denote the transpose of a matrix $M$.
The group $\gpGL{r}{\setZ}$ acts on the ring $\rgFS{\setR}{x_1,\ldots,x_r}$ of the formal power series in $r$ indeterminates over $\setR$, 
	as follows.
For any $M\in\gpGL{r}{\setZ}$ and $f\in\rgFS{\setR}{x_1,\ldots,x_r}$,
	\envHLineDefPd[1_PL_DefActionPoly]
	{
		\racFA{f}{M}{x_1,\ldots,x_r}
	}
	{
		\Fc{f}{(x_1,\ldots,x_r) \mopTP{M}}
	}
This is a right action
	(i.e., $\racF[n]{f}{M_1M_2} =  \racF{\racFrr{f}{M_1}}{M_2}$),
	and  
	extends to an action of 
	the group ring $\setZ[\gpGL{r}{\setZ}]$ of $\gpGL{r}{\setZ}$ over $\setZ$
	in a natural way  by
	\envMPd
	{
		\racFT[n]{f}{\SmN a_j M_j}
	}
	{
		\SmN a_j \racFrrT{f}{M_j}
	}
We will need to consider actions of matrices of determinant $0$ to prove our main result,
	so we adopt the non-standard definition:
	we use the transpose $\mopTP{M}$ instead of the inverse $\mopI{M}$. 
Note that
	the action of any permutation $\sig$ on $f$ is same in both definitions,
	that is,
	\envHLineCm
	{
		\racFA{f}{\sig}{x_1,\ldots,x_r}
	}
	{
		\Fc{f}{x_{\mopI{\sig}(1)},\ldots,x_{\mopI{\sig}(r)}}
	}
	because $\mopTP{\sig}=\mopI{\sig}$ in $\gpGL{r}{\setZ}$.
	
Let $\geP$, $\geQ$, $\geR$, $\geS$ be the matrices in $\gpGL{4}{\setZ}$ given by
	\envPLineCm[1_PL_DefMatPQRS]
	{
		\geP
		=
		\MatF{1&1&1&1}{0&1&1&1}{0&0&1&1}{0&0&0&1}
	,\;
		\geQ
	=
		\MatF{1&0&1&1}{0&1&1&1}{0&0&1&1}{0&0&0&1}
	,\;
		\geR
	=
		\MatF{1&0&1&0}{0&1&1&0}{0&0&1&0}{0&0&0&1}
	,\;
		\geS
	=
		\MatF{1&0&0&1}{0&1&0&1}{0&0&1&1}{0&0&0&1}	
	}
	and let $\gePh$, $\gePs$ be the elements  in $\setZ[\gpGL{4}{\setZ}]$ given by
	\envHLineCFNmPd[1_PL_DefMatPhPs]
	{
		\gePh
	}
	{
		\matI + \MatF{1&0&0&0}{0&1&0&0}{0&1&-1&1}{0&0&0&1} + \MatF{1&0&0&0}{1&0&-1&1}{0&1&-1&1}{0&0&0&1} 
	}
	{\rule{0pt}{35pt}
		\gePs
	}
	{
		\matI + \MatF{1&0&0&0}{0&1&0&0}{0&0&1&0}{0&0&1&-1} 
		+ 
		\MatF{1&0&0&0}{0&1&0&0}{0&1&0&-1}{0&0&1&-1} + \MatF{1&0&0&0}{1&0&0&-1}{0&1&0&-1}{0&0&1&-1} 
	}
For any subset $H$ of $\gpSym{4}$,
	we define an element $\mpS{H}$ in $\setZ[\gpGL{4}{\setZ}]$ as 
	\envHLineDefPd
	{
		\mpS{H}
	}
	{
		\Sm[t]{\sig\in H} \sig
	}
For example,
	\envHLineCm
	{
		\mpS{\gpC{4}}
	}
	{
		\matI 
		+ 
		\MatF{0&0&0&1}{1&0&0&0}{0&1&0&0}{0&0&1&0} 
		+ 
		\MatF{0&0&1&0}{0&0&0&1}{1&0&0&0}{0&1&0&0}	
		+
		\MatF{0&1&0&0}{0&0&1&0}{0&0&0&1}{1&0&0&0}	
	}
	where $\gpC{4}=\SetO{\gpu, (1234), (13)(24), (1432)}$ is the cyclic subgroup in $\gpSym{4}$ generated by $(1234)$. 

Our main result is stated as follows.
\begin{theorem}\label{1_Thm}
We have 
	\envHLineCm[1_Thm_EqPSumFml]
	{
		\racF{ \ggfcQZ[]{} }{ \geOm }
	}
	{
		\ggfcSqZ[]{}
	}
	where 
	\envHLinePd[1_Thm_DefOm]
	{
		\geOm
	}
	{
		\geP \mpS{\gpSym{4}} - \bkR{ \gePs \gePh \geQ - \gePh \geR - \gePs \geS + \matI} \mpS{\gpC{4}}
	}	
\end{theorem}

Taking the homogeneous parts on both sides of \refEq{1_Thm_EqPSumFml} and substituting appropriate values for the variables $x_1,x_2,x_3$, and $x_4$,
	we can obtain the following formulas:
\begin{corollary}\label{1_Cor}
For any integer $l>4$,
	we have
	\envHLineCm[1_Cor_Eq1]
	{
		\tpSm{ \bfL_4 : \bfADM}{ \wgt{\bfL_4}=l } \zt{\bfL_4}
	}
	{
		\zt{l}
	}
	\vPack[10]\envHLineCm[1_Cor_Eq2]
	{
		\tpSm{ \bfL_4 : \bfADM}{ \wgt{\bfL_4}=l } \bkR[b]{ 2^{l_1+l_2+l_3-2} + 2^{l_1+l_2-2} + 2^{l_1-1} - 2^{l_2+l_3-1} - 2^{l_2-1} } \zt{\bfL_4}
	}
	{
		l\zt{l}
	}
	\vPack[10]\envHLineCm[1_Cor_Eq3]
	{
		\tpSm{ \bfL_4 : \bfADM}{ \wgt{\bfL_4}=l } \bkR[b]{ 2^{l_1} + 2^{l_3+1} }  \zt{\bfL_4}
	}
	{
		(l+3)\zt{l}
	}
	\vPack[10]\envHLinePd[1_Cor_Eq4]
	{\hspace{-10pt}
		\tpSm{ \bfL_4 : \bfADM}{ \wgt{\bfL_4}=l } \bkR[b]{ 3^{l_2}2^{l_1-1} - 3^{l_2} - 1 } 2^{l_1+l_3} \zt{\bfL_4}
	}
	{
		\bkR[g]{ \frc{12}{(l+7)(l+2)(l-3)} + 2 } \zt{l}
	}
\end{corollary}

Formula \refEq{1_Cor_Eq1} is,
	needless to say, 
	the sum formula \refEq{1_PL_EqSumFml} for QZVs.
Formulas \refEq{1_Cor_Eq2} and \refEq{1_Cor_Eq3} are the weighted sum formulas for QZVs 
	that were proved in \cite[Theorem 1.1]{GX09} and \cite[main theorem]{OEL13},
	respectively.
These facts guarantee that, for QZVs, \refEq{1_Thm_EqPSumFml} is a natural generalization of \refEq{1_PL_EqSumFml}. 
It appears that
	\refEq{1_Cor_Eq4} is new,
	and we note that in its coefficients, it has not only powers of $2$ but also powers of $3$.

Let $\ggfcZ{r}{x_1,\ldots,x_r}$ be the generating function of MZVs of depth $r$:
	\envHLineDefPd[1_PL_DefGGfcMZV]
	{
		\ggfcZ{r}{x_1,\ldots,x_r}
	}
	{
		\Sm{ \bfL_r : \bfADM } \zt{\bfL_r} x_1^{l_1-1} \cdots x_r^{l_r-1}
	}	
Let $\ggfcSZ[]{}$, $\ggfcDZ[]{}$, and $\ggfcTZ[]{}$
	denote
	$\ggfcZ[]{1}{}$, $\ggfcZ[]{2}{}$, and $\ggfcZ[]{3}{}$,
	respectively.
Note that $\ggfcQZ[]{}=\ggfcZ[]{4}{}$.
Let $\ztS{\bfL_r}$ be regularized shuffle-type multiple zeta values (RMZVs),
	which were introduced in \cite{IKZ06}. 
RMZVs are  MZVs if $\bfL_r$ are admissible,
	but 
	RMZVs are defined for non-admissible index sets,
	unlike MZVs
	(see \cite{IKZ06} for details). 
We thus define the generating function of RMZVs of depth $r$ by 	
	\envHLineCmDef[1_PL_DefGGfcRMZV]
	{
		\ggfcZS{r}{x_1,\ldots,x_r}
	}
	{
		\Sm{ \bfL_r } \ztS{\bfL_r} x_1^{l_1-1} \cdots x_r^{l_r-1}
	} 
	where the summation includes not only admissible index sets but also non-admissible  ones.
Let $\ggfcSZS[]{}$, $\ggfcDZS[]{}$, $\ggfcTZS[]{}$, and $\ggfcQZS[]{}$
	denote $\ggfcZS[]{1}{}$, $\ggfcZS[]{2}{}$, $\ggfcZS[]{3}{}$, and $\ggfcZS[]{4}{}$,
	respectively.

For a square matrix $M$ of order $4$ and an integer $k\in\SetO{1,2,3,4}$,
	we define a row operation $\opR[]{k}{}$ as follows:
	$\opR{k}{M}$ is the matrix produced from $M$ 
	by multiplying $k$-th row by $-1$ and then adding $(k-1)$-th and $(k+1)$-th rows to $k$-th row,
	where $0$-th and $5$-th rows mean the zero row vector.
Then,
	$\gePh$ and $\gePs$ are expressed as  
	\envHLineCFCmNm[1_PL_EqMatPhPs]
	{
		\gePh
	}
	{
		\matI + \opR{3}{\matI} + \opR[]{2}{} \circ \opR{3}{\matI}
	}
	{
		\gePs
	}
	{
		\matI + \opR{4}{\matI} + \opR[]{3}{} \circ \opR{4}{\matI} + \opR[]{2}{} \circ \opR[]{3}{} \circ \opR{4}{\matI}
	}
	and
	any two adjacent matrices are transitive by some $\opR[]{k}{}$.
(Note that $\opR[]{k}{} \circ \opR[]{k}{}$ is the identity operation.)
These expressions will help us to show the equations in $\setZ[\gpGL{4}{\setZ}]$ through the present paper.
It may be worth noting that
	the matrices $\opR{4}{\matI}$, $\opR[]{3}{} \circ \opR{4}{\matI}$, and $\opR[]{2}{} \circ \opR[]{3}{} \circ \opR{4}{\matI}$ in $\gePs$
	can be obtained from the matrices $\matI$, $\opR{3}{\matI}$, and $\opR[]{2}{} \circ \opR{3}{\matI}$ in $\gePh$,
	respectively,
	by adding $4$-th column to $3$-th column and then multiplying $4$-th column by $-1$.

The present paper is organized as follows.
\refSec[s]{sectTwo} and \ref{sectThree} each have two subsections.
In \refSec{sectTwoOne}, 
	we give some identities for $\ggfcSZS[]{}$, $\ggfcDZS[]{}$, $\ggfcTZS[]{}$, and $\ggfcQZS[]{}$,	
	and in \refSec{sectTwoTwo}, we discuss 
	a relation between $\ggfcQZS[]{}$ and $\ggfcQZ[]{}$.
	
We prove \refThm{1_Thm} in \refSec{sectThreeOne}, 
	by using the results obtained in \refSec{sectTwo} and the identity for RMZVs that was proved in \cite[\refThmA{1.1}]{Machide17}.
We derive \refCor{1_Cor} from \refThm{1_Thm} in \refSec{sectThreeTwo}.
We attach an appendix at the end of the paper, 
	in which
	the identities for $\ggfcDZS[]{}$ and $\ggfcTZS[]{}$ proved in \cite{GKZ06} and \cite{Machide13b}
	are 
	restated in terms of the actions of $\setZ[\gpGL{2}{\setZ}]$ and $\setZ[\gpGL{3}{\setZ}]$,
	respectively:
	they are lower depth versions of \refEq{1_Thm_EqPSumFml}.

\begin{remark}\label{1_Rem1}
\refThm{1_Thm} and \refCor{1_Cor} in the present paper are expansions of Theorems {1.1} and {1.2} in \cite{Machide12Ax}, 
	respectively;
	the results that will be stated  below  are expansions of the results following \refSectA{2.1} in \cite{Machide12Ax}.
The results of \refSectA{2.1} have been amplified in \cite{Machide17}.
\end{remark}

\section{Preliminaries} \label{sectTwo}
\subsection{Identities for $\ggfcSZS[]{}$, $\ggfcDZS[]{}$, $\ggfcTZS[]{}$, and $\ggfcQZS[]{}$} \label{sectTwoOne}
For functions $f_{r_1},\ldots,f_{r_j}$ such that each $f_{r_i}$ has  $r_i$ variables,
	we define a function of $r=r_1+\cdots+r_j$ variables $\bfX_r = (x_1, \ldots, x_r)$ by
	\envLineDefPd
	{
		\Fc{f_{r_1}\otimes f_{r_2}\otimes\cdots\otimes f_{r_j}}{\bfX_r} 
	}
	{
		\Fc{f_{r_1}}{x_1,\ldots,x_{r_1}}
		\Fc{f_{r_2}}{x_{r_1+1},\ldots,x_{r_1+r_2}}
		\cdots
		\Fc{f_{r_j}}{x_{r_1+r_2+\cdots+r_{j-1}+1},\ldots,x_r}
	} 
For example,
	\envPLinePd
	{
		\Fc{ \ggfcSZS[]{} \otimes \ggfcSZS[]{} }{\bfX_2}			\lnP{=}	\ggfcSZS{x_1} \ggfcSZS{x_2}
	\quad\text{and}\quad
		\Fc{ \ggfcDZS[]{}\otimes\ggfcSZS[]{} }{\bfX_3}	\lnP{=}	\ggfcDZS{x_1,x_2} \ggfcSZS{x_3}
	}
Let $\Gp{\sig}=\SetT{\sig^n}{n\in\setZ}$ denote the cyclic subgroup generated by a permutation $\sig$.

The purpose of this subsection is to prove the following identities. 
\begin{proposition}\label{2.1_Prop1}
We have
	\envHLineCENmePd
	{\label{2.1_Prop1_Eq1}
		\ggfcTZS[]{} \otimes \ggfcSZS[]{}
	}
	{
		\racF{ \ggfcQZS[]{} }{ \gePs \geS }
	}
	{\label{2.1_Prop1_Eq2}
		\ggfcDZS[]{} \otimes \ggfcDZS[]{} 
	}
	{
		\racF{ \ggfcQZS[]{} }{ \gePh \geR \mpS{\Gp{(13)(24)}} }
	}
	{\label{2.1_Prop1_Eq3}
		\ggfcDZS[]{} \otimes \ggfcSZS[]{} \otimes \ggfcSZS[]{}
	}
	{
		\racF{ \ggfcQZS[]{} }{ \gePs \gePh \geQ }
	}
	{\label{2.1_Prop1_Eq4}
		\ggfcSZS[]{} \otimes \ggfcSZS[]{} \otimes \ggfcSZS[]{} \otimes \ggfcSZS[]{}
	}
	{
		\racF{ \ggfcQZS[]{} }{ \geP \mpS{\gpSym{4}}  }
	}
\end{proposition}

We will need \refLem[s]{2.1_Lem1} and \ref{2.1_Lem2} below to prove \refProp{2.1_Prop1}.
We prepare some notation to state the lemmas,
	by referring to \cite{IKZ06}.

Let $\geP[r]$ be the matrix in $\gpGL{r}{\setZ}$ given by
	\envHLineCFLqqPd[2.2_PL_DefP]
	{
		\geP[r]
	}
	{
		\MatF{1&1& \cdots & 1}{ &1& \cdots &1}{ & &\roMTx{-7}{\ddots}& \vdots }{& & &1}
	}
	{
		\gePt[r]
	}
	{
		\MatF{1&& &}{ 1&\,1&&}{ \vdots&\vdots&\roMTx{-7}{\ddots}&  }{1&\,1& \cdots&1}
	}
For a function $f$ of $r$ variables,
	we define 
	\envMCmDef{
		\popSP{f}
	}{
		\racF{f}{\geP[r]}
	}
	that is,
	\envHLinePd[2.2_PL_DefOpSha]
	{
		\popSPT{f}{x_1,x_2,\ldots,x_r}
	}
	{
		\Fc{f}{ x_1+x_2+\cdots+x_r, x_2+\cdots+x_r,\ldots,x_r }
	}
We note 
	that 
	$\geP[4]$ is $\geP$ in \refEq{1_PL_DefMatPQRS},
	and
	$\popSP{f}=f$ if $r=1$.
Let $\lfgSHiT{j}{r}$ be the shuffle elements in $\setZ[\gpSym{r}]$ 
 	defined by
	\envHLineCmDef[2.1_PL_DefActShu]
	{
		\lfgSHiT{j}{r}
	}
	{
		\tpSm{\sig\in\gpSym{r}}{\sig(1)<\cdots<\sig(j) \atop \sig(j+1)<\cdots<\sig(r)} \sig
	}
	where $j$ and $r$ are integers with $1\leq j \leq r-1$.
We set 
	\envM
	{
		\lfgSHi{j}
	}
	{
		\lfgSHiT{j}{4}
	}
	for brevity.

\refLem[s]{2.1_Lem1} and \ref{2.1_Lem2} are stated as follows.
 
\begin{lemma}\label{2.1_Lem1}
The following identities hold. 
	\envHLineCECmNme
	{\label{2.1_Lem1_Eq1}
		\popSP{ \ggfcTZS[]{} } \otimes \ggfcSZS[]{} 
	}
	{
		\racF{ \popSP{ \ggfcQZS[]{} } }{ \lfgSHi{3} }
	}
	{\label{2.1_Lem1_Eq2}
		\popSP{ \ggfcDZS[]{} } \otimes \popSP{ \ggfcDZS[]{} } 
	}
	{
		\racF{ \popSP{ \ggfcQZS[]{} } }{ \lfgSHi{2} }
	}
	{\label{2.1_Lem1_Eq3}
		\popSP{ \ggfcDZS[]{} } \otimes \ggfcSZS[]{} \otimes \ggfcSZS[]{}
	}
	{
		\racF{ \popSP{ \ggfcQZS[]{} } }{ \lfgSHi{2} \mpS{\Gp{(34)}} }
	}
	{\label{2.1_Lem1_Eq4}
		\ggfcSZS[]{} \otimes \ggfcSZS[]{} \otimes \ggfcSZS[]{} \otimes \ggfcSZS[]{}
	}
	{
		\racF{ \popSP{ \ggfcQZS[]{} } }{ \mpS{\gpSym{4}} }
	}
	where the concrete expressions of 
	$\lfgSHi{3}$, $\lfgSHi{2}$, $\lfgSHi{2} \mpS{\Gp{(34)}}$ 
	are 
	\envHLineCSNmePd 
	{\label{2.1_Lem1_EqHelpA}
		\lfgSHi{3}
	}
	{
		\gpu+(34)+(234)+(1234)
	}
	{\label{2.1_Lem1_EqHelpB}
		\lfgSHi{2}
	}
	{
		\gpu+(23)+(13)(24)+(123)+(243)+(1243)
	}
	{\label{2.1_Lem1_EqHelpC}
		\lfgSHi{2} \mpS{\Gp{(34)}}
	}
	{
		\gpu+(23)+(24)+(34)+(13)(24)+(123)+(124)
		\lnAH
		+
		(234)+(243)+(1234)+(1243)+(1324)
	}
\end{lemma}
\begin{lemma}\label{2.1_Lem2}
The following equations in $\setZ[\gpGL{4}{\setZ}]$ hold.
	\envHLineCSCmNme
	{\label{2.1_Lem2_Eq1}
		\geP\,\lfgSHi{3} 
	}
	{
		\gePs \geP 
	}
	{\label{2.1_Lem2_Eq2}
		\geP\,\lfgSHiT{2}{3}  
	}
	{
		\gePh \geP 
	}
	{\label{2.1_Lem2_Eq3}
		\geP\,\lfgSHi{2} \mpS{\Gp{(34)}} 
	}
	{
		\gePs \gePh \geP 
	}
	where $\lfgSHiT{2}{3}$ in \refEq{2.1_Lem2_Eq2} is regarded as an element in $\setZ[\gpGL{4}{\setZ}]$
	by identifying $\gpSym{3}$ as the subgroup of $\gpSym{4}$ consisting of elements that fix $4$,
	that is,
	\envHLinePd[2.1_Lem1_EqHelp]
	{
		\lfgSHiT{2}{3}
	}
	{
		\gpu + (23) + (123)
	}
\end{lemma}

Let $\rgMat{r}{\setZ}$ be the ring of square matrices of order $r$ over $\setZ$.
For any pair $(A,B)$ in $\rgMat{i}{\setZ} \times \rgMat{j}{\setZ}$, 
	we define a block diagonal matrix $A\oplus B$ in  $\rgMat{i+j}{\setZ}$ by
	\envHLineDefPd
	{
		A\oplus B
	}
	{
		\MatT{A&}{&B}
	}
Since
	$\mopTP[n]{A\oplus B}=\mopTP{A}\oplus\mopTP{B}$,
	we see from the definitions of $\otimes$ and $\oplus$ that
	\envHLineCm
	{
		\racF{f \otimes g}{A\oplus B}
	}
	{
		\racF{f}{A} \otimes \racF{g}{B}
	}
	where $f$ and $g$ are functions of $i$ and $j$ variables,
	respectively.

We now prove \refProp{2.1_Prop1}.
We will then discuss proofs of \refLem[s]{2.1_Lem1} and \ref{2.1_Lem2}.

\envProof[\refProp{2.1_Prop1}]{
By definition, 
	\envOTLineTh 
	{
		f
	}
	{
		\racF{f}{ \geP[r]\gePi[r] }
	}
	{
		\racF{\popSP{f}}{\gePi[r]}
	}
	for any function $f$ of $r$ variables,
	and so 	
	\envHLineCSNmePd
	{\label{2.1_Prop1Pr_Eq1helpA}
		\ggfcTZS[]{} \otimes \ggfcSZS[]{}
	}
	{
		\racF{ \popSP{ \ggfcTZS[]{} } \otimes \ggfcSZS[]{} }{ \gePi[3] \oplus \matI[1] }
	}
	{\label{2.1_Prop1Pr_Eq2helpA}
		\ggfcDZS[]{} \otimes \ggfcDZS[]{}
	}
	{
		\racF{ \popSP{ \ggfcDZS[]{} } \otimes \popSP{ \ggfcDZS[]{} }  }{ \gePi[2] \oplus \gePi[2] }
	}
	{\label{2.1_Prop1Pr_Eq3helpA}
		\ggfcDZS[]{} \otimes \ggfcSZS[]{} \otimes \ggfcSZS[]{}
	}
	{
		\racF{ \popSP{ \ggfcDZS[]{} } \otimes \ggfcSZS[]{} \otimes \ggfcSZS[]{} }{ \gePi[2] \oplus \matI[2] }
	}
Thus we may show that
	\envHLineCSNmePd
	{\label{2.1_Prop1Pr_Eq1helpB}
		\racF{ \popSP{ \ggfcTZS[]{} } \otimes \ggfcSZS[]{} }{ \gePi[3] \oplus \matI[1] }
	}
	{
		\racF{ \ggfcQZS[]{} }{ \gePs \geS }
	}
	{\label{2.1_Prop1Pr_Eq2helpB}
		\racF{ \popSP{ \ggfcDZS[]{} } \otimes \popSP{ \ggfcDZS[]{} }  }{ \gePi[2] \oplus \gePi[2] }
	}
	{
		\racF{ \ggfcQZS[]{} }{ \gePh \geR \mpS{\Gp{(13)(24)}} }
	}
	{\label{2.1_Prop1Pr_Eq3helpB}
		\racF{ \popSP{ \ggfcDZS[]{} } \otimes \ggfcSZS[]{} \otimes \ggfcSZS[]{} }{ \gePi[2] \oplus \matI[2] }
	}
	{
		\racF{ \ggfcQZS[]{} }{ \gePs \gePh \geQ }
	}
In fact,
	equating \refEq{2.1_Prop1Pr_Eq1helpA} and \refEq{2.1_Prop1Pr_Eq1helpB} proves \refEq{2.1_Prop1_Eq1},
	equating \refEq{2.1_Prop1Pr_Eq2helpA} and \refEq{2.1_Prop1Pr_Eq2helpB} proves \refEq{2.1_Prop1_Eq2},
	and
	so on.
We easily see that \refEq{2.1_Lem1_Eq4} yields \refEq{2.1_Prop1_Eq4} since $\popSP{\ggfcQZS[]{}}=\racF{\ggfcQZS[]{}}{\geP}$.

By the definition \refEq{2.2_PL_DefP},
	we have
	\envPLineCm
	{
		\gePi[2]	\lnP{=}	\MatT{1&-1}{0&1}
	,\quad
		\gePi[3]	\lnP{=}	\MatTh{1&-1&0}{0&1&-1}{0&0&1}
	,\quad 
		\gePi[4]	\lnP{=}	\MatF{1&-1&0&0}{0&1&-1&0}{0&0&1&-1}{0&0&0&1}
	}
	and 
	\envHLineCSNmePd
	{\label{2.1_Prop1Pr_Eq1helpZ}
		\geP\,(\gePi[3] \oplus \matI[1])
	}
	{
		\geS
	}
	{\label{2.1_Prop1Pr_Eq2helpZ}
		\geP\,(\gePi[2] \oplus \gePi[2])
	}
	{
		\geR
	}
	{\label{2.1_Prop1Pr_Eq3helpZ}
		\geP\,(\gePi[2] \oplus \matI[2])
	}
	{
		\geQ
	}
We see from \refEq{2.1_Lem1_Eq1} that
	\envOTLineThCm
	{
		\racF{ \popSP{ \ggfcTZS[]{} } \otimes \ggfcSZS[]{} }{ \gePi[3] \oplus \matI[1] }
	}
	{
		\racF{ \popSP{ \ggfcQZS[]{} } }{ \lfgSHi{3} (\gePi[3] \oplus \matI[1]) }
	}
	{
		\racF{ \ggfcQZS[]{} }{ \geP\,\lfgSHi{3} (\gePi[3] \oplus \matI[1]) }
	}
	which, together with \refEq{2.1_Lem2_Eq1} and \refEq{2.1_Prop1Pr_Eq1helpZ}, verifies \refEq{2.1_Prop1Pr_Eq1helpB}.
In a similar way,
	\refEq{2.1_Prop1Pr_Eq3helpB} is obtained from
	\refEq{2.1_Lem1_Eq3}, \refEq{2.1_Lem2_Eq3}, and \refEq{2.1_Prop1Pr_Eq3helpZ}. 
It follows from \refEq{2.1_Lem1_EqHelpB} and \refEq{2.1_Lem1_EqHelp} that
	\envHLinePd[2.1_Prop1Pr_Eq3A]
	{
		\lfgSHi{2}
	}
	{
		\lfgSHiT{2}{3} \mpS{\Gp{(13)(24)}}
	}
Since 
	\envMCm
	{
		(13)(24)
	}
	{
		\MatT{&\matI[2]}{\matI[2]&}
	}
	we have
	\envPLine
	{
		(13)(24)\MatT{A&}{&B}	\lnP{=}	\MatT{&B}{A&}
	\qquad\text{and}\qquad
		\MatT{A&}{&B} (13)(24)	\lnP{=}	\MatT{&A}{B&}
	}
	for matrices $A$ and $B$ in $\rgMat{2}{\setZ}$, 
	and so
	\envHLinePd[2.1_Prop1Pr_Eq3B]
	{
		(13)(24) (\gePi[2] \oplus \gePi[2])
	}
	{
		(\gePi[2] \oplus \gePi[2]) (13)(24)
	}
Thus,
	\envHLineCmPt{\osTx{=}{\refEq{2.1_Lem1_Eq2}}}
	{
		\racF{ \popSP{ \ggfcDZS[]{} } \otimes \popSP{ \ggfcDZS[]{} }  }{ \gePi[2] \oplus \gePi[2] }
	}
	{
		\racF{ \popSP{ \ggfcQZS[]{} } }{ \lfgSHi{2} (\gePi[2] \oplus \gePi[2]) }
	\lnAHP{\osTx{=}{(\ref{2.1_Prop1Pr_Eq3A})}}
		\racF{ \popSP{ \ggfcQZS[]{} } }{ \lfgSHiT{2}{3} \mpS{\Gp{(13)(24)}} (\gePi[2] \oplus \gePi[2]) }
	\lnAHP{\osTx{=}{(\ref{2.1_Prop1Pr_Eq3B})}}
		\racF{ \popSP{ \ggfcQZS[]{} } }{ \lfgSHiT{2}{3} (\gePi[2] \oplus \gePi[2]) \mpS{\Gp{(13)(24)}}  }
	}
	which, 
	together with \refEq{2.1_Lem2_Eq2} and \refEq{2.1_Prop1Pr_Eq2helpZ}, 
	shows \refEq{2.1_Prop1Pr_Eq2helpB}. 
}

We prove \refLem{2.1_Lem1}.

\envProof[\refLem{2.1_Lem1}]{
The following identity was given in the proof of \cite[\refThmA{6}]{IKZ06} (see also \refRem{2.1_Rem1} below). 
For integers $j$ and $r$ with $1\leq j \leq r-1$,
	\envHLinePd[2.1_Lem99_EqGGfcRMZVbyIKZ]
	{
		\popSPr{\ggfcZS[]{j}{}}\otimes\popSPr{\ggfcZS[]{r-j}{}} 
	}
	{
		\racF{\popSPr{\ggfcZS[]{r}{}}}{\lfgSHiT{j}{r}}
	}

Identities \refEq{2.1_Lem1_Eq1} and \refEq{2.1_Lem1_Eq2} immediately follow 
	from \refEq{2.1_Lem99_EqGGfcRMZVbyIKZ} with $(j,r)=(3,4)$ and $(2,4)$,
	respectively. 
We see from \refEq{2.1_Lem99_EqGGfcRMZVbyIKZ} with $(j,r)=(1,2)$ that
	\envOTLineCm
	{
		\popSP{ \ggfcSZS[]{} } \otimes \popSP{ \ggfcSZS[]{} } (x_3,x_4)
	}
	{
		\popSP{ \ggfcDZS[]{} }(x_3,x_4) + \popSP{ \ggfcDZS[]{} }(x_4,x_3)
	}
	{
		\racFrr{ \popSP{ \ggfcDZS[]{} } }{ \mpS{\Gp{(34)}} } (x_3,x_4)
	}
	and so we obtain
	\envOTLinePd[2.1_Lem1Pr_Eq1]
	{
		\popSP{ \ggfcDZS[]{} } \otimes \popSP{ \ggfcSZS[]{} } \otimes \popSP{ \ggfcSZS[]{} }
	}
	{
		\popSP{ \ggfcDZS[]{} } \otimes \bkR{ \racF{ \popSP{ \ggfcDZS[]{} } }{  \mpS{\Gp{(34)}} } }
	}
	{
		\racF{  \popSP{ \ggfcDZS[]{} } \otimes \popSP{ \ggfcDZS[]{} }  }{  \mpS{\Gp{(34)}} }
	}
Identity \refEq{2.1_Lem1_Eq3} is verified by combining \refEq{2.1_Lem1_Eq2} and \refEq{2.1_Lem1Pr_Eq1}.
By \refEq{2.1_Lem99_EqGGfcRMZVbyIKZ} with $(j,r)=(1,2)$,
	we have
	\envHLineCm[2.1_Lem1Pr_Eq2]
	{
		\ggfcSZS[]{} \otimes \ggfcSZS[]{} \otimes \ggfcSZS[]{} \otimes \ggfcSZS[]{}
	}
	{
		\racF{ \popSP{ \ggfcDZS[]{} } \otimes \popSP{ \ggfcSZS[]{} } \otimes \popSP{ \ggfcSZS[]{} } }{ \mpS{\Gp{(12)}} }
	}
	and
	combining \refEq{2.1_Lem1_Eq3} and \refEq{2.1_Lem1Pr_Eq2} we obtain
	\envHLinePd
	{
		\ggfcSZS[]{} \otimes \ggfcSZS[]{} \otimes \ggfcSZS[]{} \otimes \ggfcSZS[]{}
	}
	{
		\racF{ \popSP{ \ggfcQZS[]{} } }{ \lfgSHi{2} \mpS{\Gp{(34)}}\mpS{\Gp{(12)}} }
	}	
Identity \refEq{2.1_Lem1_Eq4} holds 
	since
	we see from \refEq{2.1_Lem1_EqHelpC} that 
	\envHLineCm
	{
		\lfgSHi{2} \mpS{\Gp{(34)}} \mpS{\Gp{(12)}}
	}
	{
		\mpS{\gpSym{4}} 
	}
	and we complete the proof.
}

\begin{remark}\label{2.1_Rem1}
The notation 
	$\ipFc[]{F}{r}{\sh}{}$ in \cite{IKZ06} is equivalent to $\ggfcZS[]{r}{}$ in the present paper.
The proof of \refEq{2.1_Lem99_EqGGfcRMZVbyIKZ} in \cite{IKZ06}
	is summarized  in the proof of \cite[\refLemA{3.1}]{Machide14} (see \cite[\refEqA{3.11}]{Machide14}),
	where
	note that
	the inverse $\mopI{M}$ are used in \cite{IKZ06,Machide14} for the action of $\gpGL{r}{\setZ}$, 
	unlike the definition \refEq{1_PL_DefActionPoly}.
\end{remark}

Let $\sig$ and $M=(m_{ij})$ be elements in $\gpSym{4}$ and $\gpGL{4}{\setZ}$, 
	respectively.
By the definition of the embedding of $\gpSym{4}$ into $\gpGL{4}{\setZ}$,
	the multiplications of $\sig$ from left and right act on the rows and columns of $M$, 
	respectively,
	as follows.
	\envPLinePd[2.1_PL_EqMatPermutation]
	{
		\sig M	\lnP{=}	(m_{\sig^{-1}(i)j})
	\qquad\text{and}\qquad
		M \sig	\lnP{=}	(m_{i\sig(j)})
	} 
That is,
	$\sig M$ and $M \sig$ are the matrices produced from $M$ 
	by replacing every $i$-th row and  $j$-th column with $\sig^{-1}(i)$-th row and $\sig(j)$-th column,
	respectively.
Equation \refEq{2.1_PL_EqMatPermutation}  will be used repeatedly below.

We are now in a position to prove \refLem{2.1_Lem2}.

\envProof[\refLem{2.1_Lem2}]{
We see from \refEq{2.1_PL_EqMatPermutation} that
	\envPLineCm
	{
		\geP (34)		\lnP{=}	\MatF{1&1&1&1}{0&1&1&1}{0&0&1&1}{0&0&1&0}
	,\quad
		\geP (234)		\lnP{=}	\MatF{1&1&1&1}{0&1&1&1}{0&1&1&0}{0&0&1&0}
	,\quad
		\geP (1234)	\lnP{=}	\MatF{1&1&1&1}{1&1&1&0}{0&1&1&0}{0&0&1&0}
	}
	and so\,\footnote{It may be worth noting that $(34)(24)=(234)$ and $(34)(24)(14)=(234)(14)=(1234)$.}
	\envHLineCSCm
	{
		\geP (34)	
	}
	{
		\opR{4}{\geP}
	}
	{
		\geP (234)	
	}
	{
		\opR{3}{\geP (34)}		\lnP{=}	\opR[]{3}{}\circ \opR{4}{ \geP }
	}
	{
		\geP (1234)
	}
	{
		\opR{2}{\geP (234)}		\lnP{=}	\opR[]{2}{}\circ \opR[]{3}{}\circ \opR{4}{ \geP }
	}
	which,
	together with \refEq{1_PL_EqMatPhPs} and \refEq{2.1_Lem1_EqHelpA},
	proves \refEq{2.1_Lem2_Eq1}.
We can show \refEq{2.1_Lem2_Eq2} in a similar way,
	and 
	we omit the proof.
We have by \refEq{2.1_Lem2_Eq1} and \refEq{2.1_Lem2_Eq2}
	\envHLinePd
	{
		\geP\,\lfgSHi{3}\, \lfgSHiT{2}{3} 
	}
	{
		\gePs\geP\, \lfgSHiT{2}{3} 
	\lnP{=}
		\gePs \gePh \geP
	}
Since
	\envMCm
	{
		\lfgSHi{3}\, \lfgSHiT{2}{3} 
	}
	{
		\lfgSHi{2} \mpS{\Gp{(34)}}
	}
	we obtain \refEq{2.1_Lem2_Eq3}.
}

\subsection{A relation between $\ggfcQZS[]{}$ and $\ggfcQZ[]{}$} \label{sectTwoTwo}	
It holds that $\ggfcQZS[]{}\neq\ggfcQZ[]{}$,
	because RMZVs for non-admissible index sets are not zero in general.
For example,
	\envMPd
	{
		\ztS{1,1,1,2}
	}
	{
		-4 \zt{2,1,1,1}
	}
(We can make more examples from \cite[(5.2)]{IKZ06}.)

The purpose of this subsection is to prove \refProp{2.2_Prop1}.
\begin{proposition}\label{2.2_Prop1}
We have 
	\envHLinePd[2.2_Prop1_Eq]
	{
		\racF{ \ggfcQZS[]{}  }{ \geOm }
	}
	{
		\racF{ \ggfcQZ[]{} }{ \geOm }
	}
\end{proposition}

\begin{remark}\label{2.2_Rem1} %
Actual values of RMZVs do not affect \refEq{2.2_Prop1_Eq} as we can see in its proof.
That is,
	\refEq{2.2_Prop1_Eq} holds if we replace RMZVs with formal variables.
\end{remark}

We define a block diagonal matrix $\geJ$ in $\rgMat{4}{\setZ}$ composed of $(0)$ and $\matI[3]$ by
	\envOTLineDefPd
	{
		\geJ
	}
	{
		(0) \oplus \matI[3]
	}
	{
		\MatF{0&0&0&0}{0&1&0&0}{0&0&1&0}{0&0&0&1} 
	}
We also define four block diagonal matrices as 
	\envHLineCmDef[2.2_PL_DefMatT]
	{\qquad
		\geT{i}
	}
	{
		(0) \oplus \geT{i}'
		\qquad
		(i=1,2,3,4)
	}
	where 
	\envPLinePd[2.2_PL_DefMatT]
	{
		\geT{1}'	\lnP{=}	\MatTh{1&1&1}{0&1&1}{0&0&1}
	,\quad
		\geT{2}'	\lnP{=}	\MatTh{1&0&1}{0&1&1}{0&0&1}
	,\quad
		\geT{3}'	\lnP{=}	\MatTh{1&0&1}{0&1&1}{0&1&0}
	,\quad
		\geT{4}'	\lnP{=}	\MatTh{1&0&1}{1&0&0}{0&1&0}
	}

We require \refLem{2.2_Lem1} to prove \refProp{2.2_Prop1}.

\begin{lemma}\label{2.2_Lem1}
We have the following equations in $\setZ[\rgMat{4}{\setZ}]$:
	\envHLineCENmePd
	{\label{2.2_Lem1_Eq1}
		\geJ \geP \mpS{\gpSym{4}} 
	}
	{
		\geT{1} \mpS{\gpSym{4}}
	}
	{\label{2.2_Lem1_Eq2}
		\geJ \gePs \gePh \geQ \mpS{\gpC{4}} 
	}
	{
		\geT{1} \mpS{\gpSym{4}}  + (\geT{2} + \geT{3} + \geT{4}) (\gpu+(243))  \mpS{\gpC{4}} 
	}
	{\label{2.2_Lem1_Eq3}
		\geJ \gePh \geR \mpS{\gpC{4}} 
	}
	{
		(\geT{2} + \geT{3} + \geT{4}) (243) \mpS{\gpC{4}}
	}
	{\label{2.2_Lem1_Eq4}
		\geJ \gePs \geS \mpS{\gpC{4}}
	}
	{
		(\geT{2} + \geT{3} + \geT{4} + \geJ) \mpS{\gpC{4}} 
	}
\end{lemma}

We now prove \refProp{2.2_Prop1}.
We will then prove \refLem{2.2_Lem1}.

\envProof[\refProp{2.2_Prop1}]{
Let $\ggfcQZSj[]{}$ denoted the difference between $\ggfcQZS[]{}$ and $\ggfcQZ[]{}$,
	that is,
	\envHLineDefPd
	{
		\ggfcQZSj[]{}
	}
	{
		\ggfcQZS[]{} - \ggfcQZ[]{} 
	}
Since
	an index set $\bfL=(l_1,l_2,l_3,l_4)$ is not admissible if and only if $l_1=1$,
	we see from \refEq{1_PL_DefGGfcMZV} and \refEq{1_PL_DefGGfcRMZV} that
	\envOTLineCm
	{
		\ggfcQZS{\bfX} - \ggfcQZ{\bfX}
	}
	{
		\ggfcQZS{0,x_2,x_3,x_4}
	}
	{
		\ggfcQZS{\bfX\mopTP{\geJ}}
	}
	or
	\envHLineCm[2.2_Prop1Pr_EqGgfcQs]
	{
		\ggfcQZSj[]{}
	}
	{
		\racF{ \ggfcQZS[]{} }{ \geJ }
	}
	which,  
	together with \refEq{2.2_Lem1_Eq1}, \refEq{2.2_Lem1_Eq2}, \refEq{2.2_Lem1_Eq3}, and \refEq{2.2_Lem1_Eq4},
	yield
	\envHLineCECm
	{
		\racF{ \ggfcQZSj[]{} }{ \geP \mpS{\gpSym{4}} }
	}
	{
		\racF{ \ggfcQZS[]{} }{ \geT{1} \mpS{\gpSym{4}} }
	}
	{
		\racF{ \ggfcQZSj[]{} }{ \gePs \gePh \geQ \mpS{\gpC{4}} }
	}
	{
		\racF[n]{ \ggfcQZS[]{} }{ \geT{1} \mpS{\gpSym{4}}  + (\geT{2} + \geT{3} + \geT{4}) (\gpu+(243))  \mpS{\gpC{4}}  }
	}
	{
		\racF{ \ggfcQZSj[]{} }{ \gePh \geR \mpS{\gpC{4}} }
	}
	{
		\racF[n]{ \ggfcQZS[]{} }{ (\geT{2} + \geT{3} + \geT{4}) (243) \mpS{\gpC{4}} }
	}
	{
		\racF{ \ggfcQZSj[]{} }{ \gePs \geS \mpS{\gpC{4}} }
	}
	{
		\racF[n]{ \ggfcQZS[]{} }{ (\geT{2} + \geT{3} + \geT{4} + \geJ) \mpS{\gpC{4}}  }
	}
	respectively.
Thus,
	\envHLineThCm
	{
		\racF{ \ggfcQZSj[]{} }{ \mpS{\gpC{4}}   }
	}
	{
		\racF{ \ggfcQZS[]{} }{ \geJ \mpS{\gpC{4}}   }
	}
	{
		\racF{ \ggfcQZSj[]{} }{ \geP \mpS{\gpSym{4}} }
		-
		\racF{ \ggfcQZSj[]{} }{ \gePs \gePh \geQ \mpS{\gpC{4}} }
		+
		\racF{ \ggfcQZSj[]{} }{ \gePh \geR \mpS{\gpC{4}} }
		+
		\racF{ \ggfcQZSj[]{} }{ \gePs \geS \mpS{\gpC{4}} }
	}
	and
	\envHLineCm
	{
		\racF{ \ggfcQZSj[]{} }{ \geOm }{  }
	}
	{
		0
	}
	which
	verifies \refEq{2.2_Prop1_Eq}.	
}

\envProof[\refLem{2.2_Lem1}]{
We have
	\envOFLineFCm[2.2_Lem1Pr_Eq1]
	{
		\geJ \geP 
	}
	{
		\MatF{0&0&0&0}{0&1&0&0}{0&0&1&0}{0&0&0&1}  \MatF{1&0&0&0}{1&1&0&0}{1&1&1&0}{1&1&1&1}  
	}
	{
		 \MatF{0&0&0&0}{0&1&0&0}{0&1&1&0}{0&1&1&1}  
	}
	{
		\geT{1}
	}
	which proves \refEq{2.2_Lem1_Eq1}.

We will next show \refEq{2.2_Lem1_Eq3} and \refEq{2.2_Lem1_Eq4}.
Direct calculations similar to \refEq{2.2_Lem1Pr_Eq1} yield
	\envHLineCSPd[2.2_Lem1Pr_Eq2]
	{
		\geJ \geR
	}
	{
		\geT{4} (243)
	}
	{
		\geJ \opR{3}{\geR}  
	}
	{
		\geT{3} (243)
	}
	{
		\geJ \opR[]{2}{} \circ \opR{3}{\geR} 
	}
	{
		\geT{2} (123)	
	}
Since
	\envHLineCm
	{
		(123)
	}
	{
		(243)(13)(24)
	}
	we see from \refEq{1_PL_EqMatPhPs} and \refEq{2.2_Lem1Pr_Eq2} that
	\envHLinePd[2.2_Lem1Pr_Eq3]
	{
		\geJ \gePh \geR
	}
	{
		\geT{2} (243)(13)(24) + (\geT{3} + \geT{4}) (243)
	}
Multiplying both sides of \refEq{2.2_Lem1Pr_Eq3} by $\mpS{\gpC{4}}$ from the right proves \refEq{2.2_Lem1_Eq3},
	because 
	$(13)(24) \in \gpC{4}$ and 
	\envHLinePd
	{
		(13)(24)\mpS{\gpC{4}}
	}
	{
		\mpS{\gpC{4}}
	}
In a similar way to \refEq{2.2_Lem1Pr_Eq3},
	we can obtain 
	\envHLineCm
	{
		\geJ \gePs \geS
	}
	{
		\geT{2} + \geT{3} + \geT{4} + \geJ (1234)
	}
	which,
	together with $(1234) \in \gpC{4}$,
	gives \refEq{2.2_Lem1_Eq4}.

By a straightforward calculation with \refEq{1_PL_EqMatPhPs},
	we have
	\envHLinePd
	{
		\gePs \gePh \geQ 
	}
	{
		\MatF{1&0&1&1}{0&1&1&1}{0&0&1&1}{0&0&0&1} + \MatF{1&0&1&1}{0&1&1&1}{0&0&1&1}{0&0&1&0}
		+
		\MatF{1&0&1&1}{0&1&1&1}{0&1&1&0}{0&0&1&0} + \MatF{1&0&1&1}{1&0&1&0}{0&1&1&0}{0&0&1&0}
		\lnAH
		+
		\MatF{1&0&1&1}{0&1&1&1}{0&1&0&1}{0&0&0&1} + \MatF{1&0&1&1}{0&1&1&1}{0&1&0&1}{0&1&0&0}
		+
		\MatF{1&0&1&1}{0&1&1&1}{0&1&1&0}{0&1&0&0} + \MatF{1&0&1&1}{1&0&1&0}{0&1&1&0}{0&1&0&0}
		\lnAH
		+ 
		\MatF{1&0&1&1}{1&0&0&1}{0&1&0&1}{0&0&0&1} + \MatF{1&0&1&1}{1&0&0&1}{0&1&0&1}{0&1&0&0}
		+
		\MatF{1&0&1&1}{1&0&0&1}{1&0&0&0}{0&1&0&0} + \MatF{1&0&1&1}{1&0&1&0}{1&0&0&0}{0&1&0&0} 
	}
Thus,
	\envHLineThCm[2.2_Lem1Pr_Eq4]
	{
		\geJ \gePs \gePh \geQ 
	}
	{
		\geT{1}  + \geT{1}(34) + \geT{1}(234) + \geT{2}(1234)
		\lnAH
		+
		\geT{1}(23)  + \geT{1}(243) + \geT{1}(24) + \geT{3}(1234) 
		\lnAH
		+ 
		\geT{2}(123)  + \geT{3}(123)  + \geT{4}(123) + \geT{4}(1234) 
	}
	{
		\geT{1} \mpS{\gpST}  + \bkR{ \geT{2}+\geT{3}+\geT{4} }  \bkR{ (1234) + (123) }
	}
	where $\gpST$ is the permutation group on the set $\Set{2,3,4}$.
Multiplying the first and last lines of \refEq{2.2_Lem1Pr_Eq4} by $\mpS{\gpC{4}}$ from the right,
	we obtain \refEq{2.2_Lem1_Eq2},
	and 
	complete the proof.
}

\section{Proofs} \label{sectThree}
\subsection{Proof of \refThm{1_Thm}} \label{sectThreeOne}
We begin by showing \refLem{3.1_Lem1}.
\begin{lemma}\label{3.1_Lem1}
Let $\gpCs{4}$ be the subset $\Set{\gpu,(1234)}$ in $\gpC{4}$.
We have
	\envHLineCm[3.1_Lem1_Eq]
	{
		\ggfcSZS[]{}^{\otimes4} 
		-
		\racF{ \ggfcDZS[]{}\otimes \ggfcSZS[]{}^{\otimes2} }{ \mpS{\gpC{4}} }
		+
		\racF{ \ggfcDZS[]{}^{\otimes2} }{  \mpS{\gpCs{4}}}
		+
		\racF{ \ggfcTZS[]{}\otimes\ggfcSZS[]{} }{ \mpS{\gpC{4}} }
		-
		\racF{ \ggfcQZS[]{} }{ \mpS{\gpC{4}} }
	}
	{
		\ggfcSqZ[]{}
	}
	where 
	 $f^{\otimes r}$ means $ \usbTx{r}{ f \otimes \cdots \otimes f}$.
\end{lemma}

\envProof{
We denote by $S_\sh$, $D_\sh$, $T_\sh$, and $Q_\sh$ the functions of one, two, three, and four variables 
	given by the restrictions of $\ztS[]{}$ to the domains $\setN$, $\setN^2$, $\setN^3$, and $\setN^4$,
	respectively.	
The following identity was shown in \cite[\refThmA{1.1}]{Machide17}: 
	\envHLineCm[3.1_Lem1Pr_EqCycFmlDpt4]
	{
		S_\sh^{\otimes4} 
		-
		\racF{D_\sh\otimes S_\sh^{\otimes2}}{ \mpS{\gpC{4}} }
		+
		\racF{D_\sh^{\otimes2}}{ \mpS{\gpCs{4}} }
		+
		\racF{T_\sh\otimes S_\sh}{ \mpS{\gpC{4}} }
		-
		\racF{Q_\sh}{ \mpS{\gpC{4}} }
	}
	{
		q_\sh	
	}
	where 
	$\Fc{q_\sh}{l_1,l_2,l_3,l_4}$ is equal to $0$ if $l_1=l_2=l_3=l_4=1$ and $\zt{l_1+l_2+l_3+l_4}$ otherwise.
Multiplying both sides of \refEq{3.1_Lem1Pr_EqCycFmlDpt4} by $x_1^{l_1}x_2^{l_2}x_3^{l_3}x_4^{l_4}$ 
	and 
	summing up over all index sets $(l_1,l_2,l_3,l_4)\in\setN^4$,
	we can obtain \refEq{3.1_Lem1_Eq}.
We omit the detail since the calculation is straightforward.
}

We are now able to prove \refThm{1_Thm}.
\envProof[\refThm{1_Thm}]{
Since 
	\envOTLineCm
	{
		\mpS{\gpC{4}}
	}
	{
		(\gpu+(13)(24))(\gpu+(1234))
	}
	{
		\mpS{\Gp{(13)(24)}} \mpS{\gpCs{4}}
	}
	we see from \refProp{2.1_Prop1} and \refLem{3.1_Lem1} that
	\envHLineThCm
	{
		\racF{ \ggfcQZS[]{} }{ \geOm }
	}
	{
		\ggfcSZS[]{}^{\otimes4} 
		-
		\racF{ \ggfcDZS[]{}\otimes \ggfcSZS[]{}^{\otimes2} }{ \mpS{\gpC{4}} }
		+
		\racF{ \ggfcDZS[]{}^{\otimes2} }{  \mpS{\gpCs{4}}}
		+
		\racF{ \ggfcTZS[]{}\otimes\ggfcSZS[]{} }{ \mpS{\gpC{4}} }
		-
		\racF{ \ggfcQZS[]{} }{ \mpS{\gpC{4}} }
	}
	{
		\ggfcSqZ[]{}
	}
	which,
	together with \refProp{2.2_Prop1},
	proves \refEq{1_Thm_EqPSumFml}.
}

\subsection{Proof of \refCor{1_Cor}} \label{sectThreeTwo}
We will prove the formulas in \refCor{1_Cor}
	by taking the homogeneous parts on both sides of \refEq{1_Thm_EqPSumFml} and substituting appropriate values for the variables.
Although the method of the proof is simple,
	it requires many calculations,
	and so
	we begin by introducing some notation and terminology, which will be useful for presenting the calculations.

We begin by defining the notation and terminology that we will use to show the fomulas.
Let $\setQ[\setR^4]$ be the free module on $\setR^4$ over $\setQ$,
	where we consider elements in $\setR^4$ as row vectors such that $\bfX=(x_1,x_2,x_3,x_4)$.
For a function $\Fc{f}{\bfX}$ with the domain $\setR^4$,
	we extend it to a homomorphism with the domain $\setQ[\setR^4]$ in a natural way by
	\envHLineDefPd[2.2_PL_DefFctExt]
	{\qquad
		\Fcs{f}{\bfAA}
	}
	{
		\Sm{i} a_i \Fc{f}{\bfX_i}  
		\qquad
		(\bfAA=\Sm{i} a_i \bfX_i \in \setQ[\setR^4])
	}
For the extension,
	we assign the same symbol $f$, 
	but
	we use square brackets $\bkS{\,}$ instead of parentheses $\bkR{\,}$.
The difference between $\Fcs{f}{\bfAA}$ and $\Fc{f}{\bfAA}$ is demonstrated as follows:
	\envPLineCm
	{
		\Fcs{f}{\bfU - 2\bfV} 	\lnP{=} 	\Fc{f}{\bfU} - 2 \Fc{f}{\bfV}
	\qquad\text{and}\qquad
		\Fc{f}{\bfU - 2\bfV} 	\lnP{=}	\Fc{f}{\bfW}
	}
	where $\bfU=(1,1,1,1)$, $\bfV=(0,1,0,1)$, and $\bfW=(1,-1,1,-1)$.
The matrix ring $\rgMat{4}{\setZ}$ acts on $\setQ[\setR^4]$ by the right multiplication,
	and so 
	it acts on the extended function \refEq{2.2_PL_DefFctExt} as
	\envHLinePd[2.2_PL_DefActFctExt]
	{
		\racFAs{f}{M}{\bfAA}
	}
	{
		\Fcs{f}{\bfAA\,\mopTP{M}}
	}
This action is a generalization of \refEq{1_PL_DefActionPoly}
	since $\racFAs{f}{M}{\bfAA} = \racFA{f}{M}{\bfAA}$ if $\bfAA \in \setR^4$,
	and
	we will use the same notation,
	$\racF{f}{M}$.
For convenience,
	we extend \refEq{2.2_PL_DefActFctExt} to an action of the free module $\setQ[\rgMat{4}{\setZ}]$ in the usual way by
	\envMPd
	{
		\racFT[n]{f}{ \SmN b_j M_j }
	}
	{
		\SmN b_j \racFrrT{f}{M_j}
	}

For any integer $l>4$,
	let $\gfcQZ[]{l}{}$ and $\gfcSqZ[]{l}{}$ be the homogeneous parts of degree $l-4$ of $\ggfcQZ[]{}$ and $\ggfcSqZ[]{}$,
	respectively.
Equivalently,
	these are defined as
	\envHLineCFDefNmePd
	{\label{3.2_PL_DefGfcQZ}
		\gfcQZ{l}{\bfX_4}
	}
	{
		\tpSm{ \bfL_4 : \bfADM }{\wgt{\bfL_4}=l} \zt{l_1,l_2,l_3,l_4} x_1^{l_1-1} x_2^{l_2-1} x_3^{l_3-1} x_4^{l_4-1} 
	}
	{\label{3.2_PL_DefFcQS}
		\gfcSqZ{l}{\bfX_4}
	}
	{
		\zt{l} \tpSm{\bfL_4}{ \wgt{\bfL_4}=l } x_1^{l_1-1} x_2^{l_2-1} x_3^{l_3-1} x_4^{l_4-1}
	}
We obtain from \refEq{1_Thm_EqPSumFml} that
	\envHLinePd[3.2_PL_EqPSumFml]
	{
		\racF{ \gfcQZ[]{l}{} }{ \geOm }
	}	
	{
		\gfcSqZ[]{l}{}
	}	
Let $\mdE=\mdE_\setQ$ be the submodule in $\setQ[\setZ^4]$ defined by 
	\envHLineDefPd
	{
		\mdE
	}
	{
		\bkB[B]{ \SmN r_i\bfK_i \mVert[b] \text{$r_i \in \setQ$ and $\bfK_i = (0,b_i,c_i,d_i) \in \setZ^4$} }
	}
We define a congruence relation $\equiv$ 
	such that 
	$\bfAA\equiv\bfBB$ if and only if $\bfAA-\bfBB\in\mdE$, 
	where $\bfAA,\bfBB\in\setQ[\setZ^4]$.
We have
	\envHLine
	{
		\gfcQZ[]{l}{\bkS{\bfAA}}
	}
	{
		\gfcQZ[]{l}{\bkS{\bfBB}}
	}
	when $\bfAA\equiv\bfBB$,
	since
	$\gfcQZ{l}{0,b,c,d} = 0$ for integers $b,c,d \in \setZ$.
For short,
	a row vector $(a,b,c,d) \in \setZ^4$ will be expressed as
	\envHLineDefPd
	{
		\veE{abcd}
	}
	{
		(a,b,c,d)
	}

We require \refProp{3.2_Prop1} to prove \refCor{1_Cor}.

\begin{proposition}\label{3.2_Prop1}
We have the following congruence equations in $\setQ[\setZ^4]$:
	\envHLineCmPt[3.2_Prop1_Eq1]{\equiv}
	{
		\veE{1000} \, \mopTP{ \geOm }
	}
	{
		\veE{1111}
	}
	\vPack[10]\envHLineCmPt[3.2_Prop1_Eq2]{\equiv}
	{
		\veE{1100} \, \mopTP{ \geOm }
	}
	{
		2\bkR{ \veE{2221}-\veE{1221} } + \bkR{ \veE{2211}+\veE{2111}-\veE{1211} } - 3\veE{1111} 
	}
	\vPack[10]\envHLineCmPt[3.2_Prop1_Eq3]{\equiv}
	{
		\bkR[B]{ \veE{1100}-\frc{2}{1}\veE{1010} } \, \mopTP{ \geOm }
	}
	{
		\veE{2111} + 2\veE{1121} -3\veE{1111} 
	}
	\vPack[10]\envHLinePdPt[3.2_Prop1_Eq4]{\equiv}
	{
		\bkR[B]{ \veE{1100}+\frc{2}{1}\veE{1010}+\veE{1110}+\frc{4}{1}\veE{1111} } \, \mopTP{ \geOm }
	}
	{
		6\bkR{ \veE{4321}-\veE{2321} } - 2\veE{2121} -\veE{1111} 
	}
\end{proposition}

We now prove \refCor{1_Cor}.
We will then discuss a proof of \refProp{3.2_Prop1}.

\envProof[\refCor{1_Cor}]{
Let $\setN_0$ denote the set of non-negative integers,
	and
	let $\nmSet{X}$ denote the number of the elements of a set $X$.
Considering the correspondence
	\envHLinePt{\,\leftrightarrow\,}
	{
		k_1+k_2+\cdots+k_s
	}
	{
		\usbTx{k_1}{\circ\cdots\circ}\,|\,\usbTx{k_2}{\circ\cdots\circ}\,|\,\cdots\,|\,\usbTx{k_s}{\circ\cdots\circ}
	} 
	for integers $k_1, k_2, \ldots, k_s \in \setN_0$,
	we see that 
	\envHLineCm[3.2_1_CorPr_EqHelp1]
	{
		\nmSet{ \SetT{(k_1,\ldots,k_s) \in \setN_0^r}{ k_1+\cdots+k_s=k } }
	}
	{
		\binom{k+s-1}{s-1}
	}
	where $k\geq0$ and $s\geq1$.

Let $\bfX_4=(x_1,x_2,x_3,x_4)$ be a vector in $\Set{0,1}^4$ with $\bfX_4\neq(0,0,0,0)$,
	and
	let $i_1, \ldots, i_r$ be the distinct indices such that $x_{i_1}=\cdots=x_{i_r}=0$:
	note that $0 \leq r \leq 3$.
We have
	\envLineFCmPte
	{
		\tpSm{\bfL_4}{ \wgt{\bfL_4}=l } x_1^{l_1-1} x_2^{l_2-1} x_3^{l_3-1} x_4^{l_4-1}
	}{ = }
	{
		\nmSet{ \SetT{\bfL_4  \in \setN^4}{ \text{$\wgt{\bfL_4}=l$ and $l_{i_1}=\cdots=l_{i_r}=1$}} } 
	}{ \osTx{=}{(l \leftrightarrow k+1)} }
	{
		\nmSet{ \SetT{\bfK_4 \in \setN_0^4}{ \text{$\wgt{\bfK_4}=l-4$ and $k_{i_1}=\cdots=k_{i_r}=0$}} }  
	}{ = }
	{
		\nmSet{ \SetT{\bfK_{4-r} \in \setN_0^{4-r}}{ \text{$\wgt{\bfK_{4-r}}=l-4$}} }
	}
	which, together with \refEq{3.2_1_CorPr_EqHelp1}, gives
	\envHLinePd
	{
		\tpSm{\bfL_4}{ \wgt{\bfL_4}=l } x_1^{l_1-1} x_2^{l_2-1} x_3^{l_3-1} x_4^{l_4-1}
	}
	{
		\binom{l-r-1}{3-r}
	}
By the definition \refEq{3.2_PL_DefFcQS},
	we thus obtain 
	\envHLine[3.2_1_CorPr_Eq0help]
	{
		\gfcSqZ{l}{\bfX}
	}
	{
		\envCaseFPd[d]{
			\zt{l}
			&
			(\bfX=\veE{1000})
		}{
			(l-3)\zt{l}
			&
			(\bfX=\veE{1100}, \veE{1010})
		}{\rule{0pt}{20pt}
			\frc{2}{(l-2)(l-3)}\zt{l}
			&
			(\bfX=\veE{1110})
		}{\rule{0pt}{20pt}
			\frc{6}{(l-1)(l-2)(l-3)}\zt{l}
			&
			(\bfX=\veE{1111})
		}
	}	

We now prove the desired formulas.
We see from \refEq{3.2_Prop1_Eq1} and \refEq{3.2_1_CorPr_Eq0help} that
	\envHLineCFLaaCm
	{
		\gfcQZ[]{l}{ \bkS{ \veE{1000} \, \mopTP{ \geOm } } }
	}
	{
		\gfcQZ{l}{\veE{1111}}
	}
	{
		\gfcSqZ{l}{\veE{1000}}
	}
	{
		\zt{l}
	}
	respectively,
	which, together with \refEq{3.2_PL_EqPSumFml}, 
	yields
	\envHLinePd[3.2_1_CorPr_Eqq1]
	{
		\gfcQZ{l}{\veE{1111}}
	}
	{
		\zt{l}
	}
Rewriting the left-hand side of \refEq{3.2_1_CorPr_Eqq1},
	we obtain \refEq{1_Cor_Eq1}.
	
In a similar way,
	it follows from \refEq{3.2_Prop1_Eq2} and \refEq{3.2_1_CorPr_Eq0help} that
	\envLine
	{
		\gfcQZ[]{l}{ \bkS{ \veE{1100} \, \mopTP{ \geOm } } }
	}
	{
		2\bkR{ \gfcQZ{l}{\veE{2221}}-\gfcQZ{l}{\veE{1221}} } + \gfcQZ{l}{\veE{2211}}+\gfcQZ{l}{\veE{2111}}-\gfcQZ{l}{\veE{1211}} - 3 \gfcQZ{l}{\veE{1111}}
	}
	and
	\envHLineCm
	{
		\gfcSqZ{l}{\veE{1100}}
	}
	{
		(l-3)\zt{l}
	}
	respectively,
	and so
	\envHLineCm[3.2_1_CorPr_Eqq2]
	{
		2\bkR{ \gfcQZ{l}{\veE{2221}}-\gfcQZ{l}{\veE{1221}} } + \gfcQZ{l}{\veE{2211}}+\gfcQZ{l}{\veE{2111}}-\gfcQZ{l}{\veE{1211}}  
	}
	{
		l\zt{l}
	}
	which shows \refEq{1_Cor_Eq2}. 
	
Combining \refEq{3.2_Prop1_Eq3} and \refEq{3.2_1_CorPr_Eq0help},
 	together with \refEq{3.2_PL_EqPSumFml},
 	yields 
	\envHLinePd[3.2_1_CorPr_Eqq3]
	{
		\gfcQZ{l}{\veE{2111}} + 2\gfcQZ{l}{\veE{1121}} 
	}
	{
		\frc{2}{l+3}\zt{l}
	}
Multiplying both sides of \refEq{3.2_1_CorPr_Eqq3} by $2$, 
	we prove \refEq{1_Cor_Eq3}.

We can deduce from \refEq{3.2_Prop1_Eq4} and \refEq{3.2_1_CorPr_Eq0help} that
	\envLinePd[3.2_1_CorPr_Eqq4help]
	{
		6\bkR{ \gfcQZ{l}{\veE{4321}}-\gfcQZ{l}{\veE{2321}} } - 2\gfcQZ{l}{\veE{2121}} 
	}
	{
		\bkR[g]{ (l-3) + \frc{2}{l-3} + \frc{2}{(l-2)(l-3)} + \frc{24}{(l-1)(l-2)(l-3)} + 1} \zt{l}
	}
We see that
	\envHLineThCm
	{
		\text{(LHS of \refEq{3.2_1_CorPr_Eqq4help})}
	}
	{
		\tpSm{ \bfL_4 : \bfADM}{ \wgt{\bfL_4}=l } \bkR[b]{ 4^{l_1-1}3^{l_2}2^{l_3} - 3^{l_2}2^{l_1+l_3-1} - 2^{l_1+l_3-1} }  \zt{\bfL} 
	}
	{
		\tpSm{ \bfL_4 : \bfADM}{ \wgt{\bfL_4}=l } \bkR[b]{ 3^{l_2}2^{l_1-1} - 3^{l_2} - 1 } 2^{l_1+l_3-1} \zt{\bfL}
	}
	and
	\envHLineThPd
	{
		\text{(RHS of \refEq{3.2_1_CorPr_Eqq4help})}
	}
	{
		\bkR[g]{ \frc{24}{ \bkB{ 36 + 12(l-2) + (l-1)(l-2) }(l-3)} + 1} \zt{l}
	}
	{
		\bkR[g]{ \frc{24}{(l+7)(l+2)(l-3)} + 1 }\zt{l}
	}
Thus,
	multiplying both sides of \refEq{3.2_1_CorPr_Eqq4help} by $2$,
	we obtain \refEq{1_Cor_Eq4},
	which
	completes the proof.
}

We prepare \refLem{3.2_Lem1} to prove \refProp{3.2_Prop1}.

\begin{lemma}\label{3.2_Lem1}
We have the following congruence equations in $\setQ[\setZ^4]$:
\mbox{}\\\bfTx{(i)}
	\envLinePt[3.2_Lem1_Eq1]{\equiv} 
	{
		\bfX \, \mopTP[n]{ \geP \mpS{\gpSym{4}} }
	}
	{
		\envCaseFPd{
			6\bkR{ \veE{1111} + \veE{1110} + \veE{1100} + \veE{1000} }
			&
			(\bfX=\veE{1000})
		}{
			4\bkR{ \veE{2221} + \veE{2211} + \veE{2210} + \veE{2111} + \veE{2110} + \veE{2100} }
			&
			(\bfX=\veE{1100}, \veE{1010})
		}{
			6\bkR{ \veE{3321} + \veE{3221} + \veE{3211} + \veE{3210} }
			&
			(\bfX=\veE{1110})
		}{
			24 \veE{4321}
			&
			(\bfX=\veE{1111})
		}
	}
\bfTx{(ii)}
	\envLinePt[3.2_Lem1_Eq2]{\equiv} 
	{
		\bfX \, \mopTP[n]{ \gePs \gePh \geQ \mpS{\gpC{4}} }
	}
	{
		\envCaseFiPd{
			12\veE{1000} + 10\veE{1100} + 8\veE{1110} + 6\veE{1111} 
			&
			(\bfX=\veE{1000})
		}{\rule{0pt}{15pt}
			6\bkR{ \veE{2100}+\veE{1111} } + 5\veE{2110} 
			\\
			+ 
			4\bkR{ \veE{2210}+\veE{2111}+\veE{1110} } + 3\veE{2211} 
			\\
			+ 
			2\bkR{ \veE{2221}+\veE{1221} + \veE{1211}+\veE{1210}+\veE{1121}+\veE{1100} }
			\\
			+ 
			\veE{1101} + \veE{1011} + \veE{1010} + \veE{1001} 
			&
			(\bfX=\veE{1100})
		}{\rule{0pt}{15pt}
			8\veE{2100} + 6\veE{2110} 
			\\
			+
			4\bkR{ \veE{2210}+\veE{2111}+\veE{1221} + \veE{1211}+\veE{1210}+\veE{1121} }
			\\
			2\bkR{ \veE{2211}+\veE{1101}+\veE{1011}+\veE{1010}+\veE{1001} }
			&
			(\bfX=\veE{1010})
		}{\rule{0pt}{15pt}
			8\veE{2221} + 6\bkR{ \veE{3210}+\veE{2321}+\veE{2211} } 
			\\
			+
			4\bkR{ \veE{3211}+\veE{2210}+\veE{2121}+\veE{2111} } 
			\\
			+ 
			2\bkR{ \veE{3221}+\veE{2110}+\veE{2101} } 
			&
			(\bfX=\veE{1110})
		}{\rule{0pt}{15pt}
			24\veE{3321} + 16\veE{3221} + 8\veE{3211}
			&
			(\bfX=\veE{1111})
		}
	}
\bfTx{(iii)}
	\envLinePt[3.2_Lem1_Eq3]{\equiv}
	{
		\bfX \, \mopTP[n]{ \gePh \geR \mpS{\gpC{4}} }
	}
	{
		\envCaseFiPd{
			3\veE{1000} + 2\veE{1100} + \veE{1110} 
			&
			(\bfX=\veE{1000})
		}{\rule{0pt}{15pt}
			3\veE{1111} + 2\bkR{ \veE{1210}+\veE{1110} } 
			\\
			+ 
			\veE{1211}+\veE{1100}+\veE{1011}+\veE{1010}+\veE{1001}   
			&
			(\bfX=\veE{1100})
		}{\rule{0pt}{15pt}
			4\veE{2100} + 2\veE{2110} 
			&
			(\bfX=\veE{1010})
		}{\rule{0pt}{15pt}
			2\bkR{ \veE{2210}+\veE{2111}+\veE{1221}+\veE{1121} } 
			\\
			+ 
			\veE{2211}+\veE{2110}+\veE{1211}+\veE{1101} 
			&
			(\bfX=\veE{1110})
		}{\rule{0pt}{15pt}
			8\veE{2221} + 4\veE{2211}
			&
			(\bfX=\veE{1111})
		}
	}
\bfTx{(iv)}
	\envLinePt[3.2_Lem1_Eq4]{\equiv} 
	{
		\bfX \, \mopTP[n]{ \gePs \geS \mpS{\gpC{4}} }
	}
	{
		\envCaseFiPd{
			4\veE{1000} + 2\veE{1100} + \veE{1111}+\veE{1110} 
			&
			(\bfX=\veE{1000})
		}{\rule{0pt}{15pt}
			2\bkR{ \veE{2100}+\veE{1121}+\veE{1110}+\veE{1100} } 
			\\
			+ 
			\veE{2111}+\veE{2110}+\veE{1101}+\veE{1001} 
			&
			(\bfX=\veE{1100})
		}{\rule{0pt}{15pt}
			4\bkR{ \veE{1210}+\veE{1010} } + 2\bkR{ \veE{1211}+\veE{1101}+\veE{1011}+\veE{1001} } 
			&
			(\bfX=\veE{1010})
		}{\rule{0pt}{15pt}
			3\veE{1111} + 2\bkR{ \veE{2210}+\veE{2121}+\veE{2101}+\veE{1221} }
			\\
			+
			\veE{2211}+\veE{2110}+\veE{1211}+\veE{1110}+\veE{1011} 
			&
			(\bfX=\veE{1110})
		}{\rule{0pt}{15pt}
			8\veE{2221} + 4\bkR{ \veE{2211}+\veE{2111} }
			&
			(\bfX=\veE{1111})
		}
	}
\bfTx{(v)}
	\envHLinePt[3.2_Lem1_Eq5]{\equiv} 
	{
		\bfX \, \mopTP{\mpS{\gpC{4}}} 
	}
	{
		\envCaseFiPd{
			\veE{1000} 
			&
			(\bfX=\veE{1000})
		}{
			\veE{1100} + \veE{1001}
			&
			(\bfX=\veE{1100})
		}{
			2 \veE{1010} 
			&
			(\bfX=\veE{1010})
		}{
			\veE{1110} + \veE{1101} + \veE{1011} 
			&
			(\bfX=\veE{1110})
		}{
			4\veE{1111}
			&
			(\bfX=\veE{1111})
		}
	}
\end{lemma}

Before discussing a proof of \refLem{3.2_Lem1}, 
	we will show \refProp{3.2_Prop1}
	by substituting 
	the congruence equations in \refLem{3.2_Lem1} into 
	\envHLine[3.2_PL_EqXOmega]
	{
		\bfX \, \mopTP{ \geOm }
	}
	{
		\bfX \, 
		\bkB{ 
			\mopTP[n]{ \geP \mpS{\gpSym{4}} } 
			- 
			\mopTP[n]{ \gePs \gePh \geQ \mpS{\gpC{4}} } 
			+ 
			\mopTP[n]{ \gePh \geR \mpS{\gpC{4}} }
			+
			\mopTP[n]{ \gePs \geS \mpS{\gpC{4}} }
			-
			\mopTP{\mpS{\gpC{4}}} 
		}
	}
	for $\bfX\in\Set{\veE{1000},\veE{1100},\veE{1010},\veE{1110},\veE{1111}}$.
We note that 
	the right-hand sides of 
	\refEq{3.2_Lem1_Eq1}, \refEq{3.2_Lem1_Eq2}, \refEq{3.2_Lem1_Eq3}, \refEq{3.2_Lem1_Eq4}, and \refEq{3.2_Lem1_Eq5} in \refLem{3.2_Lem1}
	include vectors in 
	\envPLine{
		\setV = \SetT{\veE{i_1i_2i_3i_4}}{i_j=0\,(\exists j)};
	}
	however, 
	those of \refEq{3.2_Prop1_Eq1}, \refEq{3.2_Prop1_Eq2}, \refEq{3.2_Prop1_Eq3}, and \refEq{3.2_Prop1_Eq4} in \refProp{3.2_Prop1}  
	do not include such vectors.
That is,
	in calculating \refEq{3.2_PL_EqXOmega}, 
	the vectors in $\setV$ cancel each other.
This fact will help us with the proof of \refProp{3.2_Prop1}.

\envProof[\refProp{3.2_Prop1}]{
Substituting equations from \refEq{3.2_Lem1_Eq1} through \refEq{3.2_Lem1_Eq5} for $\bfX=\veE{1000}$
	into the right-hand side of \refEq{3.2_PL_EqXOmega}, 
	we obtain 
	\envLineThCmPt{\equiv}
	{
		\veE{1000} \, \mopTP{ \geOm }
	}
	{
		6\bkR{ \veE{1111} + \veE{1110} + \veE{1100} + \veE{1000} } - \bkR{ 12\veE{1000} + 10\veE{1100} + 8\veE{1110} + 6\veE{1111}   }
		\lnAH
		+
		\bkR{ 3\veE{1000} + 2\veE{1100} + \veE{1110}  } + \bkR{ 4\veE{1000} + 2\veE{1100} + \veE{1111}+\veE{1110}  }
		-
		\veE{1000} 
	}
	{
		\veE{1111}
	}
	which proves \refEq{3.2_Prop1_Eq1}.
	
By \refEq{3.2_Lem1_Eq2} and \refEq{3.2_Lem1_Eq3} for $\bfX=\veE{1100}$,
	we have 
	\envLineCmPt[3.2_Prop1Pr_Eq1]{\equiv}
	{
		\veE{1100} \, \bkB{ \mopTP[n]{ \gePs \gePh \geQ \mpS{\gpC{4}} } - \mopTP[n]{ \gePh \geR \mpS{\gpC{4}} } }
	}
	{
		6\veE{2100} + 5\veE{2110} + 4\bkR{ \veE{2210}+\veE{2111} } + 3\bkR{ \veE{2211}+\veE{1111} } 
		\lnAH
		+
		2\bkR{ \veE{2221}+\veE{1221}+\veE{1121}+\veE{1110} } + \veE{1211} + \veE{1101} + \veE{1100}  
	}
	and
	by \refEq{3.2_Lem1_Eq4} and \refEq{3.2_Lem1_Eq5} for $\bfX=\veE{1100}$,
	\envLinePdPt[3.2_Prop1Pr_Eq2]{\equiv}
	{
		\veE{1100} \, \bkB{ \mopTP[n]{ \gePs \geS \mpS{\gpC{4}} } - \mopTP{\mpS{\gpC{4}}} }
	}
	{
		2\bkR{ \veE{2100}+\veE{1121}+\veE{1110} } + \veE{2111}+\veE{2110}+\veE{1101} + \veE{1100} 
	}
Noting \refEq{3.2_Prop1Pr_Eq1} and \refEq{3.2_Prop1Pr_Eq2},
	we see that 
	substituting equations from \refEq{3.2_Lem1_Eq1} through \refEq{3.2_Lem1_Eq5} for $\bfX=\veE{1100}$ into the right-hand side of \refEq{3.2_PL_EqXOmega}
	gives
	\envHLineFCmPt{\equiv}
	{
		\veE{1100} \, \mopTP{ \geOm }
	}
	{
		\veE{1100} \, \mopTP[n]{ \geP \mpS{\gpSym{4}} } 
		\lnAH
		-
		\veE{1100} \, \bkB{ \mopTP[n]{ \gePs \gePh \geQ \mpS{\gpC{4}} } - \mopTP[n]{ \gePh \geR \mpS{\gpC{4}} } }
		+
		\veE{1100} \, \bkB{ \mopTP[n]{ \gePs \geS \mpS{\gpC{4}} } - \mopTP{\mpS{\gpC{4}}} }
	}
	{
		4\bkR{ \veE{2221} + \veE{2211} + \veE{2210} + \veE{2111} + \veE{2110} + \veE{2100} }
		\lnAH
		-
		6\veE{2100} - 5\veE{2110} - 4\bkR{ \veE{2210}+\veE{2111} } - 3\bkR{ \veE{2211}+\veE{1111} } 
		\lnAH
		-
		2\bkR{ \veE{2221}+\veE{1221}+\veE{1121}+\veE{1110} } - \veE{1211} - \veE{1101} - \veE{1100}  
		\lnAH
		+
		2\bkR{ \veE{2100}+\veE{1121}+\veE{1110} } + \veE{2111}+\veE{2110}+\veE{1101} + \veE{1100} 
	}
	{
		-3\veE{1111} + 2\bkR{ \veE{2221}-\veE{1221} } + \veE{2211}+\veE{2111}-\veE{1211} 
	}
	which proves \refEq{3.2_Prop1_Eq2}.

Similarly,
	we obtain
	\envLineCmPt{\equiv}
	{
		\veE{1010} \, \mopTP{ \geOm }
	}
	{
		4\bkR{ \veE{2221} + \veE{2211} + \veE{2210} + \veE{2111} + \veE{2110} + \veE{2100} }
		\lnAH
		-
		\bkB{
			8\veE{2100} + 6\veE{2110} 
			+
			4\bkR{ \veE{2210}+\veE{2111}+\veE{1221} + \veE{1211}+\veE{1210}+\veE{1121} }
			\lnAHs{10}
			+
			2\bkR{ \veE{2211}+\veE{1101}+\veE{1011}+\veE{1010}+\veE{1001} }
		}
		\lnAH[]
		+
		4\veE{2100} + 2\veE{2110} + 4\bkR{ \veE{1210}+\veE{1010} } 
		\lnAH[]
		+ 
		2\bkR{ \veE{1211}+\veE{1101}+\veE{1011}+\veE{1001} }  - 2 \veE{1010} 
	}
	which can be summarized  as  
	\envHLinePdPt[3.2_Prop1Pr_Eq3]{\equiv}
	{
		\veE{1010} \, \mopTP{ \geOm }
	}
	{
		4\bkR{ \veE{2221}-\veE{1221}-\veE{1121} } + 2\bkR{ \veE{2211}-\veE{1211} }
	}
Combining \refEq{3.2_Prop1_Eq2} and \refEq{3.2_Prop1Pr_Eq3} yields
	\envHLineThCmPt{\equiv}
	{
		\bkR[B]{ \veE{1100}-\frc{2}{1}\veE{1010} } \, \mopTP{ \geOm }
	}
	{
		2\bkR{ \veE{2221}-\veE{1221} } + \veE{2211}+\veE{2111}-\veE{1211} - 3\veE{1111} 
		\lnAH
		-
		2\bkR{ \veE{2221}-\veE{1221}-\veE{1121} } - \veE{2211} + \veE{1211} 
	}
	{
		\veE{2111} + 2\veE{1121} -3\veE{1111} 
	}
	which proves \refEq{3.2_Prop1_Eq3}.

We see from \refEq{3.2_Lem1_Eq1} and \refEq{3.2_Lem1_Eq2} for $\bfX=\veE{1110}$ that
	\envLineCmPt[3.2_Prop1Pr_Eq4]{\equiv}
	{
		\veE{1110} \, \bkB{ \mopTP[n]{ \geP \mpS{\gpSym{4}} }- \mopTP[n]{ \gePs \gePh \geQ \mpS{\gpC{4}} } }
	}
	{
		- 
		8\veE{2221} 
		+
		6 \bkR{ \veE{3321} - \veE{2321} - \veE{2211} }
		\lnAH
		+
		4 \bkR{ \veE{3221} - \veE{2210} - \veE{2121} - \veE{2111} }
		+
		2 \bkR{ \veE{3211} - \veE{2110} - \veE{2101}  }
	}
	and from \refEq{3.2_Lem1_Eq3}, \refEq{3.2_Lem1_Eq4}, and \refEq{3.2_Lem1_Eq5} for $\bfX=\veE{1110}$ that
	\envLinePdPt[3.2_Prop1Pr_Eq5]{\equiv}
	{
		\veE{1110} \, \bkB{ \mopTP[n]{ \gePh \geR \mpS{\gpC{4}} } + \mopTP[n]{ \gePs \geS \mpS{\gpC{4}} } - \mopTP{\mpS{\gpC{4}}} }
	}
	{
		4\bkR{ \veE{2210}+\veE{1221} }
		+
		3 \veE{1111}
		\lnAH
		+
		2\bkR{ \veE{2211}+\veE{2121}+\veE{2111}+\veE{2110}+\veE{2101}+\veE{1211}+\veE{1121} } 
	}
Substituting \refEq{3.2_Prop1Pr_Eq4} and \refEq{3.2_Prop1Pr_Eq5} into the right-hand side of \refEq{3.2_PL_EqXOmega}, 
	we obtain
	\envHLinePdPt[3.2_Prop1Pr_Eq6]{\equiv}
	{
		\veE{1110} \, \mopTP{ \geOm }
	}
	{
		- 8\veE{2221} +  6\bkR{ \veE{3321}-\veE{2321} } + 4\bkR{ \veE{3221}+\veE{1221}-\veE{2211} } 
		\lnAH
		+
		3\veE{1111} + 2\bkR{ \veE{3211}+\veE{1211}+\veE{1121}-\veE{2121}-\veE{2111} }
	}
Substituting equations from \refEq{3.2_Lem1_Eq1} through \refEq{3.2_Lem1_Eq5} for $\bfX=\veE{1111}$
	into the right-hand side of \refEq{3.2_PL_EqXOmega}, 
	we also have
	\envLinePdPt[3.2_Prop1Pr_Eq7]{\equiv}
	{
		\veE{1111} \, \mopTP{ \geOm }
	}
	{
		24\bkR{ \veE{4321}-\veE{3321} } + 16\bkR{ \veE{2221}-\veE{3221} } + 8\bkR{ \veE{2211}-\veE{3211} } + 4\bkR{ \veE{2111}-\veE{1111} }
	}
Combining \refEq{3.2_Prop1_Eq2}, \refEq{3.2_Prop1Pr_Eq3}, \refEq{3.2_Prop1Pr_Eq6}, and \refEq{3.2_Prop1Pr_Eq7} yields
	\envLineThCmPt{\equiv}
	{
		\bkR[B]{ \veE{1100}+\frc{2}{1}\veE{1010}+\veE{1110}+\frc{4}{1}\veE{1111} } \, \mopTP{ \geOm }
	}
	{
		\bkB{ 2\bkR{ \veE{2221}-\veE{1221} } + \bkR{ \veE{2211}+\veE{2111}-\veE{1211} } - 3\veE{1111} } 
		\lnAH
		+
		\bkB{ 2\bkR{ \veE{2221}-\veE{1221}-\veE{1121} } + \veE{2211}-\veE{1211} }
		\lnAH
		+
		\bkB{
			- 
			8\veE{2221} +  6\bkR{ \veE{3321}-\veE{2321} } + 4\bkR{ \veE{3221}+\veE{1221}-\veE{2211} } 
                		\lnAHs{10}
                		+
                		3\veE{1111} + 2\bkR{ \veE{3211}+\veE{1211}+\veE{1121}-\veE{2121}-\veE{2111} }
		}
		\lnAH
		+
		\bkB{ 6\bkR{ \veE{4321}-\veE{3321} } + 4\bkR{ \veE{2221}-\veE{3221} } + 2\bkR{ \veE{2211}-\veE{3211} } + \veE{2111}-\veE{1111} }
	}
	{
		6\bkR{ \veE{4321}-\veE{2321} } - 2\veE{2121} -\veE{1111} 
	}
	which proves \refEq{3.2_Prop1_Eq4},
	and this completes the proof.
}

For a matrix $M=(m_{ij})$ in $\rgMat{4}{\setZ}$ and distinct integers $i_1,\ldots,i_r$ in $\Set{1,2,3,4}$ ($r\leq4$),
	we define a row vector in $\setZ^4$ by 
	\envHLineCmDef
	{
		\VMp{M}{i_1\ldots i_r}
	}
	{
		\veE{\delta_1\delta_2\delta_3\delta_4}M
	}
	where $\delta_i$ is $1$ if $i\in\Set{i_1,\ldots,i_r}$ and $0$ otherwise.
Equivalently,
	$\VMp{M}{i_1\ldots i_r}$ is determined by
	\envHLinePd
	{
		\VMp{M}{i_1\ldots i_r}
	}
	{
		\bkR[g]{ \SmT{k=1}{r} m_{i_k1}, \SmT{k=1}{r} m_{i_k2}, \SmT{k=1}{r} m_{i_k3}, \SmT{k=1}{r} m_{i_k4} }
	}
For example,
	if 
	\envPLineCm
	{
		M	\lnP{=}	\gePt		\lnP{=}	\MatF{1&0&0&0}{1&1&0&0}{1&1&1&0}{1&1&1&1}
	}
	then
	\envPLinePd{
		\VMp{M}{2}	\lnP{=}	\veE{0100}M	\lnP{=}	\veE{1100}
	\qquad\text{and}\qquad
		\VMp{M}{124}	\lnP{=}	\veE{1101}M	\lnP{=}	\veE{3211}
	}

Recall that a permutation $\sig$ in $\gpSym{4}$ is identified with the matrix $(\delta_{i\sig(j)})_{1\leq i,j\leq 4}$.
Thus,
	\envMCm
	{
		\bfX \sig 
	}
	{
		(x_{\sig(1)},x_{\sig(2)},x_{\sig(3)},x_{\sig(4)})
	}
	and
	we can deduce the following equations by direct calculations:
	\envHLine[3.2_Lem2Pr_Eqq1]
	{
		\bfX \, \mpS{\gpSym{4}} 
	}
	{
		\envCaseFCm{
			6\bkR{ \veE{1000} + \veE{0100} + \veE{0010} + \veE{0001} }
			&
			(\bfX=\veE{1000})
		}{
			4\bkR{ 
				\veE{1100} + \veE{1010} + \veE{1001} 
				\\ \hspace{10pt}
				+  
				\veE{0110} + \veE{0101} + \veE{0011} 
			}
			&
			(\bfX=\veE{1100}, \veE{1010})
		}{
			6\bkR{ \veE{1110} + \veE{1101} + \veE{1011} + \veE{0111} }
			&
			(\bfX=\veE{1110})
		}{
			24\veE{1111}
			&
			(\bfX=\veE{1111})
		}
	}	
	and
	\envHLine[3.2_Lem2Pr_Eqq2]
	{
		\bfX \, \mpS{\gpC{4}} 
	}
	{
		\envCaseFiPd{
			\veE{1000} + \veE{0100} + \veE{0010} + \veE{0001} 
			&
			(\bfX=\veE{1000})
		}{
			\veE{1100} + \veE{0110} + \veE{0011} + \veE{1001}
			&
			(\bfX=\veE{1100})
		}{
			2\bkR{ \veE{1010} + \veE{0101} }
			&
			(\bfX=\veE{1010})
		}{
			\veE{1110} + \veE{1101} + \veE{1011} + \veE{0111} 
			&
			(\bfX=\veE{1110})
		}{
			4\veE{1111}
			&
			(\bfX=\veE{1111})
		}
	}
For any $M\in\rgMat{4}{\setZ}$,
	we thus have
	\envHLine[3.2_Lem2_Eq1]
	{\hspace{-15pt}
		\bfX \, \mpS{\gpSym{4}} M 
	}
	{
		\envCaseFCm{
			6\bkR{ \VMp{M}{1} + \VMp{M}{2} + \VMp{M}{3} + \VMp{M}{4} }
			&
			(\bfX=\veE{1000})
		}{\rule{0pt}{15pt}
			4\bkR{ 
				\VMp{M}{12} +\VMp{M}{13} + \VMp{M}{14} 
				\\\hspace{10pt}
				+ 
				\VMp{M}{23} + \VMp{M}{24} + \VMp{M}{34} 
			}
			&
			(\bfX=\veE{1100}, \veE{1010})
		}{\rule{0pt}{15pt}
			6\bkR{ \VMp{M}{123} + \VMp{M}{124} + \VMp{M}{134} + \VMp{M}{234} }
			&
			(\bfX=\veE{1110})
		}{\rule{0pt}{15pt}
			24\VMp{M}{1234}
			&
			(\bfX=\veE{1111})
		}
	}
	and
	\envHLine[3.2_Lem2_Eq2]
	{
		\bfX \, \mpS{\gpC{4}} M 
	}
	{
		\envCaseFiPd{
			\VMp{M}{1} + \VMp{M}{2} + \VMp{M}{3} + \VMp{M}{4}
			&
			(\bfX=\veE{1000})
		}{\rule{0pt}{15pt}
			\VMp{M}{12} + \VMp{M}{23} + \VMp{M}{34} + \VMp{M}{14}
			&
			(\bfX=\veE{1100})
		}{\rule{0pt}{15pt}
			2\bkR{ \VMp{M}{13} + \VMp{M}{24}  }
			&
			(\bfX=\veE{1010})
		}{\rule{0pt}{15pt}
			\VMp{M}{123} + \VMp{M}{124} + \VMp{M}{134} + \VMp{M}{234}
			&
			(\bfX=\veE{1110})
		}{\rule{0pt}{15pt}
			4\VMp{M}{1234}
			&
			(\bfX=\veE{1111})
		}
	}

Using \refEq{3.2_Lem2_Eq1} and \refEq{3.2_Lem2_Eq2},
	we now prove \refLem{3.2_Lem1}
	for the completeness of the proof of \refProp{3.2_Prop1},
	or 
	for that of \refCor{1_Cor}.

\envProof[\refLem{3.2_Lem1}]{
We obtain by \refEq{3.2_Lem2_Eq1} 
	\envHLineTh
	{
		\bfX \, \mopTP[n]{ \geP \mpS{\gpSym{4}} }  
	}
	{
		\bfX \, \mpS{\gpSym{4}} \mopTP{ \geP}
	}
	{
		\envCaseFCm{
			6\bkR{ \veE{1000} + \veE{1100} + \veE{1110} + \veE{1111} }
			&
			(\bfX=\veE{1000})
		}{
			4\bkR{ \veE{2100} + \veE{2110} + \veE{2111} + \veE{2210} + \veE{2211} + \veE{2221} }
			&
			(\bfX=\veE{1100}, \veE{1010})
		}{
			6\bkR{ \veE{3210} + \veE{3211} + \veE{3221} + \veE{3321} }
			&
			(\bfX=\veE{1110})
		}{
			24 \veE{4321}
			&
			(\bfX=\veE{1111})
		}
	}
	which proves \refEq{3.2_Lem1_Eq1}.

We can obtain by \refEq{3.2_Lem2_Eq2} the following congruence equations:
	\envHLinePt[3.2_Lem1Pr_Eqqiii_2]{\equiv}
	{
		\bfX \, \mpS{\gpC{4}} \MatF{1&0&0&0}{0&1&0&0}{1&1&1&0}{0&0&0&1} 
	}
	{
		\envCaseFiCm{
			\veE{1000}+\veE{1110} 
			&
			(\bfX=\veE{1000})
		}{\rule{0pt}{15pt}
			\veE{1100}+\veE{1210}+\veE{1111}+\veE{1001} 
			&
			(\bfX=\veE{1100})
		}{\rule{0pt}{15pt}
			2 \veE{2110}
			&
			(\bfX=\veE{1010})
		}{\rule{0pt}{15pt}
			\veE{2210}+\veE{1101}+\veE{2111}+\veE{1211} 
			&
			(\bfX=\veE{1110})
		}{\rule{0pt}{15pt}
			4\veE{2211}
			&
			(\bfX=\veE{1111})
		}
	}
	\envHLinePt[3.2_Lem1Pr_Eqqiii_3]{\equiv}
	{
		\bfX \, \mpS{\gpC{4}} \MatF{1&0&0&0}{0&1&1&0}{1&1&0&0}{0&0&1&1} 
	}
	{
		\envCaseFiCm{
			\veE{1000}+\veE{1100} 
			&
			(\bfX=\veE{1000})
		}{\rule{0pt}{15pt}
			\veE{1110}+\veE{1210}+\veE{1111}+\veE{1011} 
			&
			(\bfX=\veE{1100})
		}{\rule{0pt}{15pt}
			2\veE{2100} 
			&
			(\bfX=\veE{1010})
		}{\rule{0pt}{15pt}
			\veE{2210}+\veE{1121}+\veE{2111}+\veE{1221} 
			&
			(\bfX=\veE{1110})
		}{\rule{0pt}{15pt}
			4\veE{2221} 
			&
			(\bfX=\veE{1111})
		}
	}
	\envHLinePt[3.2_Lem1Pr_Eqqiii_4]{\equiv}
	{
		\bfX \, \mpS{\gpC{4}} \MatF{1&1&0&0}{0&0&1&0}{1&0&0&0}{0&1&1&1} 
	}
	{
		\envCaseFiPd{
			\veE{1100}+\veE{1000}
			&
			(\bfX=\veE{1000})
		}{\rule{0pt}{15pt}
			\veE{1110}+\veE{1010}+\veE{1111}+\veE{1211}
			&
			(\bfX=\veE{1100})
		}{\rule{0pt}{15pt}
			2 \veE{2100} 
			&
			(\bfX=\veE{1010})
		}{\rule{0pt}{15pt}
			\veE{2110}+\veE{1221} +\veE{2211}+\veE{1121}
			&
			(\bfX=\veE{1110})
		}{\rule{0pt}{15pt}
			4\veE{2221} 
			&
			(\bfX=\veE{1111})
		}
	}
We see from \refEq{1_PL_EqMatPhPs} that
	\envHLineCm
	{
		\mopTP[n]{ \gePh \geR }
	}
	{
		\MatF{1&0&0&0}{0&1&0&0}{1&1&1&0}{0&0&0&1} 
		+ 
		\MatF{1&0&0&0}{0&1&1&0}{1&1&0&0}{0&0&1&1} 
		+ 
		\MatF{1&1&0&0}{0&0&1&0}{1&0&0&0}{0&1&1&1}
	}
	and so
	\envHLineThPd
	{
		\bfX \, \mopTP[n]{ \gePh \geR \mpS{\gpC{4}} } 
	}
	{
		\bfX \, \mpS{\gpC{4}} \mopTP[n]{ \gePh \geR }
	}
	{
		\bfX \, \mpS{\gpC{4}} \MatF{1&0&0&0}{0&1&0&0}{1&1&1&0}{0&0&0&1}  
		+ 
		\bfX \, \mpS{\gpC{4}} \MatF{1&0&0&0}{0&1&1&0}{1&1&0&0}{0&0&1&1}  
		+ 
		\bfX \, \mpS{\gpC{4}} \MatF{1&1&0&0}{0&0&1&0}{1&0&0&0}{0&1&1&1} 
	}
Thus, 
	the sum of \refEq{3.2_Lem1Pr_Eqqiii_2}, \refEq{3.2_Lem1Pr_Eqqiii_3}, and \refEq{3.2_Lem1Pr_Eqqiii_4}
	proves
	\refEq{3.2_Lem1_Eq3}.

With the same method that we used for \refEq{3.2_Lem1_Eq3},
	we can prove \refEq{3.2_Lem1_Eq2} and \refEq{3.2_Lem1_Eq4}.
We omit these proofs 
	because 
	of space limitations.

Since $\veE{abcd}\equiv0$ for $a=0$,		
	\refEq{3.2_Lem1_Eq5} immediately follows from \refEq{3.2_Lem2_Eq2} with $M=\matI$,
	which completes the proof.
}

\section*{Acknowledgements}
This work was supported by JST ERATO Grant Number JPMJER1201, Japan.
\mbox{}\\

\renewcommand{\thesubsubsection}{Appendix \Alph{subsubsection}}
\subsubsection{} \label{sectAppA}	
Let $\gpC{3}$ denote the cyclic subgroup $\Gp{(123)}$ in $\gpSym{3}$.
We note that $\gpC{3}$ is the alternating group $\gpAlt{3}$ of degree $3$.
We define the formal power series $\ggfcSdZ{x_1,x_2}$ and $\ggfcStZ{x_1,x_2,x_3}$
	as
	\envHLineCFCmDef
	{
		\ggfcSdZ{x_1,x_2}
	}
	{
		\tpSm{\bfL_2}{ \wgt{\bfL_2}>2} \zt{ l_1+l_2 } x_1^{l_1-1} x_2^{l_2-1} 
	}
	{
		\ggfcStZ{x_1,x_2,x_3}
	}
	{
		\tpSm{\bfL_3}{ \wgt{\bfL_3}>3} \zt{ l_1+l_2+l_3 } x_1^{l_1-1} x_2^{l_2-1} x_3^{l_3-1} 
	}
	respectively.
	
The generalizations of \refEq{1_PL_EqSumFml} for DZVs and TZVs from the viewpoint of the generating functions $\ggfcDZ{x_1,x_2}$ and $\ggfcTZ{x_1,x_2,x_3}$
	can be stated as follows: 
	\envHLineCm[App_PL_EqPSumFmlDZVhelp]
	{
		\ggfcSdZ{x_1,x_2}
	}
	{
		\Sm{\sig\in\gpSym{2}} \bkR{ \ggfcDZ{x_{\sig(1)}+x_{\sig(2)},x_{\sig(2)}} - \ggfcDZ{x_{\sig(1)},x_{\sig(2)}} }
	}
	\vPack\envHLinePd[App_PL_EqPSumFmlTZVhelp]
	{
		\ggfcStZ{x_1,x_2,x_3}
	}
	{
		\Sm{\sig\in\gpSym{3}} \ggfcTZ{x_{\sig(1)}+x_{\sig(2)}+x_{\sig(3)},x_{\sig(2)}+x_{\sig(3)},x_{\sig(3)}}
		\lnAH
		-
		\Sm{\sig\in\gpC{3}}
		\bkR[G]{
			\Sm{\tau\in\Gp{(23)}} \ggfcTZ{x_{\sig(1)}+x_{\sig(3)},x_{\sig\tau(2)}+x_{\sig\tau(3)},x_{\sig\tau(3)}}
			\lnAHs{40}
			+
			\ggfcTZ{x_{\sig(1)}+x_{\sig(2)},x_{\sig(2)},x_{\sig(3)}}
			-
			\ggfcTZ{x_{\sig(1)},x_{\sig(2)},x_{\sig(3)}}
		}
	}
The above formulas were proved in \cite[(27)]{GKZ06} and \cite[\refThmA{1.2}]{Machide13b},
	respectively.
Note that
	the original formula in \cite[\refThmA{1.2}]{Machide13b} was not written in terms of formal power series but in terms of homogeneous polynomials.

Let $\geQ[3]$, $\geR[3]$ be the matrices in $\gpGL{3}{\setZ}$ given by
	\envPLineCm
	{
		\geQ[3]
	\lnP{=}
		\MatTh{1&1&0}{0&1&0}{0&0&1}
	,\qquad
		\geR[3]
	\lnP{=}
		\MatTh{1&0&1}{0&1&1}{0&0&1}	
	}
	and let $\gePh[3]$ be the element in $\setZ[\gpGL{3}{\setZ}]$ given by
	\envHLinePd
	{
		\gePh[3]
	}
	{
		\matI[3] + \MatTh{1&0&0}{0&1&0}{0&1&-1} 
	}
Recall 
	from \refEq{2.2_PL_DefP} that
	\envPLinePd
	{
		\geP[2]	\lnP{=}	\MatT{1&1}{0&1}
	,\qquad
		\geP[3]	\lnP{=}	\MatTh{1&1&1}{0&1&1}{0&0&1}
	}
	
We can now restate \refEq{App_PL_EqPSumFmlDZVhelp} and \refEq{App_PL_EqPSumFmlTZVhelp} as follows.
	
\begin{theoremN}[\cite{GKZ06,Machide13b}]\label{App_Thm}
We have 
	\envHLineCm[App_Thm_Eq1]
	{
		\racF{ \ggfcDZ[]{} }{ \bkR{ \geP[2] - \matI[2] } \mpS{\gpSym{2}}  }
	}
	{
		\ggfcSdZ[]{}
	}
	\vPack[20]\envHLinePd[App_Thm_Eq2]
	{
		\racF[n]{ \ggfcTZ[]{} }{ \geP[3] \mpS{\gpSym{3}} - \bkR{ \geQ[3] + \gePh[3] \geR[3] - \matI[3]} \mpS{\gpC{3}} }
	}
	{
		\ggfcStZ[]{}
	}
\end{theoremN}



\end{document}